\documentclass[final,hidelinks,onefignum,onetabnum]{siamart251216}

\usepackage{amssymb}
\usepackage{graphicx}
\usepackage{float}
\usepackage{enumerate}
\usepackage{multirow}
\usepackage{booktabs}
\usepackage{array}
\usepackage{threeparttable}

\newsiamremark{remark}{Remark}
\newsiamremark{hypothesis}{Hypothesis}
\newsiamthm{conj}{Conjecture}
\newsiamremark{exmp}{Example}
\crefname{hypothesis}{Hypothesis}{Hypotheses}
\crefname{conj}{Conjecture}{Conjectures}
\crefname{exmp}{Example}{Examples}

\newcommand{\rmd}{\mathrm{d}}
\newcommand{\vx}{\mathbf{x}}

\newcommand{\R}{{\mathbb R}}

\headers{A PP method for NLSE with low-regularity potential}{W. Bao, Q. Sun, and C. Wang}
\title{A potential projection method for the nonlinear Schr\"{o}dinger equation with low regularity potential}
\author{Weizhu Bao\thanks{Department of Mathematics, National University of Singapore, Singapore 119076
(\email{matbaowz@nus.edu.sg}, \email{qingyuan.sun@u.nus.edu}).
}
\and
Qingyuan Sun\footnotemark[1]
\thanks{Zhiyuan College and School of Mathematical Sciences, Shanghai Jiao Tong University, Shanghai 200240, China.
}
\and
Chushan Wang\thanks{Department of Statistics and CCAM, University of Chicago, Chicago, IL 60637
(\email{chushanwang@uchicago.edu}).
}
}

\ifpdf
\hypersetup{
  pdftitle={A potential projection method for the nonlinear Schr\"{o}dinger equation with low regularity potential},
  pdfauthor={Anonymous authors}
}
\fi

\begin{document}
\maketitle

\begin{abstract}
We propose a potential projection (PP) method for the nonlinear Schr\"{o}dinger equation (NLSE) with low-regularity and singular potentials by replacing the potential with its Fourier projection (truncated Fourier series). This gives rise to an NLSE with projected potential (NLSEwPP). We rigorously establish its well-posedness under low-regularity potentials, and its convergence to the original NLSE. The NLSEwPP is then directly discretized by the standard time-splitting Fourier pseudospectral (TS-FP) method with the mesh size chosen according to the potential projection parameter, resulting in the fully discrete TS-PP-FP method. Numerically, TS-PP-FP consistently achieves optimal temporal and spatial convergence rates with respect to the regularity of the exact solution under low-regularity and even singular potentials. Compared with the time-splitting Fourier spectral method, TS-PP-FP attains the same convergence rates at significantly lower computational cost and preserves the mass conservation and time symmetry (that are lost in the Fourier spectral discretization of the time-splitting method); compared with the time-splitting FP method, it is substantially more accurate for low-regularity potentials with the same computational cost. 
\end{abstract}



% \begin{highlights}
% \item A Fourier pseudo-spectral method with potential projection is proposed for low-regularity and singular potentials.
% \item The projected equation converges with rates $O(N^{-2})$ in $L^2$ and $O(N^{-1})$ in $H^1$.
% \item The method avoids the enlarged auxiliary grid of eFP while retaining the efficiency and structure of FP.
% \item Numerical tests demonstrate accuracy, and efficiency for rough and singular potentials.
% \end{highlights}

% REQUIRED

% REQUIRED

\section{Introduction}
Consider the following nonlinear Schr\"{o}dinger equation (NLSE):
	\begin{equation}\label{eq::nlse} 
\begin{cases}\:i\partial_t\psi(\mathbf{x},t)=-\Delta\psi(\mathbf{x},t)+V(\mathbf{x})
\psi(\mathbf{x},t)+f(|\psi|^2)\psi(\mathbf{x},t),\quad\mathbf{x}\in\Omega,\:t>0,\\
\:\psi(\mathbf{x},0)=\psi_0(\mathbf{x}),\quad\mathbf{x}\in\overline{\Omega},\end{cases}
	\end{equation}
	where $t$ denotes time, $\mathbf{x}=(x_1,\ldots,x_d)^T \in \mathbb{R}^d$ ($d=1,2,3$) is the spatial 
coordinate, $\psi := \psi(\mathbf{x},t)$ is a complex-valued wave function, and $\Omega = \prod_{j=1}^d(a_j,b_j) \subset \mathbb{R}^d$ is a bounded periodic domain. Here, $V:=V(\mathbf{x}): \Omega \rightarrow \mathbb{R}$ is a time-independent real-valued potential, possibly of low regularity or singular, and $f(\rho) = \beta \rho$ with a given constant $\beta \in \mathbb R$. There are two important conserved quantities of the NLSE \cref{eq::nlse}: the mass and energy given by
    \begin{align}
        &M(\psi(\cdot, t)): = \int_{\Omega} |\psi(\vx, t)|^2 \rmd \vx \equiv  M(\psi(\cdot, 0)), \quad t \geq 0,  \label{eq:mass}\\
        &E(\psi(\cdot, t)):=\int_{\Omega}\left[|\nabla\psi|^2
        +V(\vx)|\psi|^2+\frac{\beta}{2}|\psi|^4\right]\rmd \vx
        \equiv E(\psi(\cdot,0)),\quad t\geq0. \label{eq:energy}
    \end{align}
    In addition, \cref{eq::nlse} is time symmetric (or time reversible) and gauge invariant (or time-transverse invariant) \cite{antoine2013nlsgpe}. 
	
	The NLSE \eqref{eq::nlse}  is a central model in quantum physics and chemistry, nonlinear optics, and materials science \cite{fibich2015nonlinear,sulem1999nonlinear,cazenave2003semilinear}. In particular, it is widely used to model and simulate Bose-Einstein condensates (BEC), where it is also known as the Gross-Pitaevskii equation (GPE) \cite{erdos2007derivation,bao2013bec,BEC2003}. In many classical settings, the potential $V$ is smooth; typical examples include the harmonic potential $V(\vx)=|\vx|^2/2$ and the optical lattice potential \cite{antoine2013nlsgpe}. 
    % which reads, in 3D, $V(\vx)=\cos(Lx) + \cos(Ly) + \cos(Lz)$, where $\vx = (x, y, z)^T$ and $L \in \R$ is a constant. 
    For such smooth problems, many accurate and efficient numerical methods have been developed. For the temporal discretizations, these include finite difference methods \cite{akrivis1993cubic,henning2017cnls,sanzserna1984methods,bao2013bec}, exponential wave integrators (EWI) \cite{celledoni2008sei,hochbruck2010expi,bao2024ewi}, time-splitting methods \cite{besse2002splitting,mclachlan2002splitting,jahnke2000error,bao2003numerical,lubich2008splitting,bao2013bec,blanes2024splitting,descombes2013lietrotter}, and, more recently, low-regularity integrators (LRI) \cite{LRI2018,LRItree,LRIreview} for low-regularity initial data. These temporal integrators are commonly combined with spatial discretizations such as finite difference methods \cite{akrivis1993cubic,sanzserna1984methods}, spectral and pseudospectral methods \cite{bao2003numerical,antoine2013nlsgpe,bao2013bec,bao2024ewi}, and finite element methods \cite{henning2017cnls,karakashian1993optimal,sanzserna1984methods}. Among these, the time-splitting Fourier pseudospectral (TS-FP) method is especially attractive as it is accurate and efficient and preserves many important dynamical properties of the continuous NLSE \cite{antoine2013nlsgpe}.
	
	Although smooth potentials are common, many physically important systems involve potentials with limited regularity or genuine singularities. Such potentials arise in Bose--Einstein condensates, quantum many-body systems, disordered media, and localized wave phenomena. One important class consists of nonsmooth but bounded potentials. Examples include the square-well potentials and step potentials, which are bounded but discontinuous and, in general, belong only to $L^\infty$ \cite{potensquarewell,fritsch2020creating,seaman2005effect}. Other examples include power-law traps $V(\vx)=|\vx|^\gamma$ with $\gamma \in (0, \infty)$, which are used to model BEC traps with different geometries and may also have only limited Sobolev regularity \cite{BECtrap,bayindir1999bose,jaouadi2011bose}. Highly irregular disordered potentials also play a central role in Anderson localization and wave propagation in random media \cite{anderson1958absence,anderson1989,whitenoisepotential}. In other fundamental models, the potential is singular. A prototypical example is the Coulomb potential $V(\vx)=-1/|\vx|$, which models the Coulomb force between charged particles in quantum mechanics and is unbounded near the origin. Further examples include inverse-square potentials, screened Coulomb potentials, and Dirac $\delta$-potentials, each of which arises in important physical models and presents additional analytical and numerical challenges \cite{nisoli2014attractive,yukawa1955interaction,zhou2023nlsdelta}.

	In the presence of low-regularity or singular potentials, although the TS-FP method can still be used, its temporal and spatial convergence orders are substantially reduced \cite{prr,bao2024efp,fang2026,becker2025Trotter}. Recently, the Fourier spectral (FS) discretization of the time-splitting methods is analyzed in \cite{bao2024tsnls}, which demonstrates significantly improved accuracy (subject to a necessary CFL-type condition) compared to the TS-FP method under low-regularity potentials. More recent works also demonstrate the potential optimal convergence of the FS method for singular potentials \cite{bao2026singular,bao2026singularoptimal}. Despite its favorable accuracy, the FS discretization has two important practical drawbacks. First, it requires accurate evaluation of the Fourier coefficients of products such as $V \phi$, where $\phi$ is the numerical solution at each time step. For low-regularity potentials, this is generally more involved and more expensive than the pointwise multiplication used in the standard FP method. Second, the structure preservation properties of the time-splitting methods (such as mass conservation and time symmetry) will be lost by the fully discrete scheme using the FS discretization in space. The first difficulty can be partially alleviated by the extended Fourier pseudospectral (eFP) method \cite{bao2024efp}, which may be understood as an efficient realization of the Fourier spectral discretization through an extended FFT. However, eFP requires an FFT on an enlarged grid at every time step, leading to a computational cost at least $3^d$ times that of the standard FP method in $d$ dimensions. This overhead becomes particularly significant in three dimensions. Moreover, eFP does not address the loss of time symmetry and mass conservation caused by the Fourier spectral treatment of the potential term.

    The main objective of this paper is to develop a new spatial discretization--the potential projection Fourier pseudospectral (PP-FP) method--for the time-splitting methods, which enjoys the advantages of the FS and FP methods at the same time. Roughly speaking, it can achieve (at least) the same accuracy as the FS method with the same computational cost and the same structure preservation (when combined with time-splitting methods) as the FP method. We also remark here that although we focus on applying the PP-FP method to time-splitting methods in this paper, it is independent of the temporal discretization and can also be easily combined with other integrators such as finite difference methods or exponential integrators. 
    
    % The PP-FP method retains the favorable convergence behavior of the FS method for low-regularity potentials. At the same time, it avoids evaluating Fourier coefficients of products such as $V\psi^n$ at each time step, and therefore retains the computational efficiency of the standard FP method. Moreover, because the projection is performed at the equation level, the resulting semidiscrete problem has the standard pseudospectral form. Consequently, when PP-FP is coupled with time-splitting methods, the fully discrete scheme inherits structural properties of the underlying splitting integrator, including mass conservation and time reversibility, which are generally lost under the FS and eFP discretizations. This feature also shows the potential in applying the time-splitting PP-FP method in the quantum computing setting. Finally, the PP-FP method is independent of the temporal integrator and can be combined with a broad range of time discretizations, including finite difference methods, time-splitting methods, and exponential-type integrators. In this work, we focus on its coupling with time-splitting methods, which are most popular in practice. 
        
    There are two key ingredients in the design of the TS-PP-FP method: (i) we first perform a potential projection that replaces the potential $V$ in the NLSE by its Fourier projection $P_N V$, which yields the NLSE with potential projection (NLSEwPP) introduced in \cref{eq::nlsem}, and (ii) we directly discretize the NLSEwPP using the standard TS-FP method with the mesh size chosen according to the projection parameter $N$. More precisely, the FP method uses $M$ degrees of freedom in each spatial direction, and in the proposed TS-PP-FP method, we take $M=N$. The relation $M\geq N$ is essential: the spatial grid must sufficiently resolve the projected potential $P_N V$ (see \cref{fig:section11}). Thus, after preprocessing to compute $P_N V$, the TS-PP-FP method has the same FFT-based implementation and the same per-step computational cost as the standard TS-FP method. Moreover, because the projection is performed at the equation level, when coupled with time-splitting methods, the PP-FP method enjoys the same structure-preservation property as the standard FP method, including mass conservation and time reversibility. This feature also shows the potential in applying the time-splitting PP-FP method in the quantum computing setting \cite{prx2021,childs2022,fang2026}. 
	
	In terms of the error estimate, 
    % we do not provide a complete error analysis for the PP-FP method under low-regularity potentials in this work. Instead, 
    in \cref{thm::error_regularization}, we prove the $H^2$ well-posedness of the NLSEwPP uniformly in the projection parameter $N$ for low-regularity potentials $V \in L^\infty$, 
    % (with $\sup_{N \geq 1} \| P_N V \|_{L^\infty} < \infty$)
    and establish the convergence of the NLSEwPP solution to the solution of the original NLSE. 
    % The convergence order is $O(N^{-2})$ in the $L^2$-norm and $O(N^{-1})$ in the $H^1$- for $L^\infty$-potential. In view of the sharp $H^2$-regularity of the exact solution under low-regularity potential $V \in L^2$ or $V \in L^\infty$, this rate is already remarkable. Moreover, for more regular potential $V \in H^\alpha$ (with $\alpha > d/2$), the convergece rates improve to $O(N^{-(2+\alpha)})$ in the $L^2$-norm and $O(N^{-(1+\alpha)})$ in the $H^1$-norm. The sharpness of these convergence rates is also confirmed by the numerical results in \cref{fig:section21}. 
    These results justify designing the TS-PP-FP method based on the NLSEwPP rather than directly on the original NLSE. Moreover, we provide comprehensive numerical experiments for potentials of varying regularity, showing that, under the same CFL-type time step restriction introduced in \cite{bao2024tsnls,bao2024efp}, the TS-PP-FP method consistently achieves the optimal spatial and temporal convergence rate with respect to the regularity of the exact solution (see \cref{sec:4}). A rigorous fully discrete error estimate for the TS-PP-FP method is out of the scope of the current paper, and will be presented separately. It should be noted that our well-posedness and convergence analysis for the NLSEwPP in fact addresses the modelling error of using the NLSEwPP to approximate the NLSE, and may be useful in much broader settings than the classical grid-based numerical methods. 

    % Finally, we combine PP-FP method with a range of time-splitting schemes, from first to fourth order, to investigate numerically how their performance depends on the regularity of the potential. In particular, we examine the optimal convergence rates attainable by splitting methods for potentials of a given regularity, as well as the minimal regularity required to achieve a prescribed convergence order. Our numerical results further suggest that imposing a CFL-type time-step restriction can substantially relax the regularity requirements and improve the performance of splitting methods when combined with the PP-FP spatial discretization. This behavior is consistent with that observed for splitting methods combined with the Fourier spectral discretization in \cite{bao2024tsnls}. 
    % % \cw{summarize the numerical findings of the time-splitting methods}
	
	The remainder of this paper is organized as follows. In \cref{sec::nlsepp}, we apply the potential projection method to obtain the NLSEwPP, and establish its well-posedness and convergence towards the original NLSE. In \cref{sec::VP}, we discretize the NLSEwPP by time-splitting Fourier pseudospectral methods to obtain fully discrete TS-PP-FP schemes. Numerical results are provided in \cref{sec:4} to illustrate the effectiveness of the proposed TS-PP-FP methods for low-regularity potentials. Finally, concluding remarks are given in \cref{sec:5}. Throughout the paper, standard
notations of Sobolev spaces and corresponding norms are adopted. The notation $A \lesssim B$ is used to represent that there exists a generic constant $C > 0$, such that $|A| \le CB$. To simplify the notation, we shall sometimes omit the spatial variable $\vx$ in the presentation when there is no confusion. 

\section{A potential projection method}\label{sec::nlsepp}

In this section, we first apply the potential projection (PP) method to the NLSE to obtain the NLSE with projected potential (NLSEwPP). Then we show the (uniform) well-posedness of the NLSEwPP, and establish its convergence to the original NLSE. The convergence order is also verified numerically, indicating its optimality. This will serve as the foundation of the TS-PP-FP method developed in \cref{sec::VP}. 

%We replace the potential $V$ by its Fourier projection onto a finite-dimensional space and consider the resulting nonlinear Schr\"odinger equation with potential projection. It provides the continuous model underlying the PP-FP method and forms the well-posedness, uniform regularity, and convergence analysis.

\subsection{NLSE with projected potential}
	We begin by introducing the NLSEwPP. 
	Let $N \in 4 \mathbb{N}^+$. Define the symmetric $d$-dimensional index set
	\begin{equation}
		\mathcal{S}_N
		=
		\left\{
		\mathbf{k}=(k_1,\ldots,k_d)\in\mathbb{Z}^d:
		-\frac{N}{2}\leq k_j \leq \frac{N}{2},\ j=1,\ldots,d
		\right\}.
	\end{equation}
	For $\mathbf{k}\in\mathcal{S}_N$, denote
	\begin{equation}
		\boldsymbol{\mu}_{\mathbf{k}}
		=
		\left(
		\frac{2\pi k_1}{b_1-a_1},
		\ldots,
		\frac{2\pi k_d}{b_d-a_d}
		\right)^T,
		\quad
		\mathbf{x}-\mathbf{a}
		=
		(x_1-a_1,\ldots,x_d-a_d)^T,
	\end{equation}
	and set
	\begin{equation}
		X_N
		=
		\text{span}
		\left\{
		e^{i\boldsymbol{\mu}_{\mathbf{k}}\cdot(\mathbf{x}-\mathbf{a})}:
		\mathbf{k}\in\mathcal{S}_N
		\right\}.
	\end{equation}
	Let $P_N : L^2(\Omega) \to X_N$ be the standard $L^2$-projection onto $X_N$ defined by \cite{STW2021}
	\begin{equation}
		\label{eq::pn}
		(P_Nu)(\mathbf{x})
		=
		\sum_{\mathbf{k}\in\mathcal{S}_N}
		\widehat{u}_{\mathbf{k}}
		e^{i\boldsymbol{\mu}_{\mathbf{k}}\cdot(\mathbf{x}-\mathbf{a})},
	\end{equation}
	where $u\in L^2(\Omega)$ and $\widehat{u}_{\mathbf{k}}$ are the Fourier coefficients of $u$ given by
	\begin{equation}
		\label{eq::fourier_coeff}
		\widehat{u}_{\mathbf{k}}
		=
		\frac{1}{|\Omega|}
		\int_{\Omega}
		u(\mathbf{x})
		e^{-i\boldsymbol{\mu}_{\mathbf{k}}\cdot(\mathbf{x}-\mathbf{a})}
		\rmd\mathbf{x},
		\quad \mathbf{k}\in\mathbb{Z}^d.
	\end{equation}
	Since $\mathcal{S}_N$ is symmetric, the projection preserves real-valuedness: if $V$ is real-valued, then $P_N V$ is real-valued.
	Define the periodic Sobolev spaces as \cite{Adams}
	\begin{equation}
		H^m_{\mathrm{per}}(\Omega):=
		\left\{
		u\in L^2(\Omega):
		\sum_{\mathbf{k}\in\mathbb{Z}^d}
		(1+|\boldsymbol{\mu}_{\mathbf{k}}|^2)^m
		|\widehat{u}_{\mathbf{k}}|^2<\infty
		\right\},\quad m\geq0.
	\end{equation}
	\begin{remark}
		For integer $m\geq1$, this Fourier-coefficient definition is equivalent to the usual boundary-condition definition \cite{Adams}
		\[
		H^m_{\mathrm{per}}(\Omega)
		=
		\left\{
		u\in H^m(\Omega):
		\partial^\alpha u|_{x_j=a_j}
		=
		\partial^\alpha u|_{x_j=b_j},\
		|\alpha|\le m-1,\ j=1,\ldots,d
		\right\}.
		\]
	\end{remark}
	Then we can introduce a projected potential $P_N V$ into the NLSE, which is well-defined for all $V \in L^2(\Omega)$. This gives the following NLSEwPP
	\begin{equation}
		\label{eq::nlsem}
		\begin{cases}
		\begin{aligned}
			i\partial_t\psi_N={}&-\Delta\psi_N
			+P_NV(\mathbf{x})\psi_N+f(|\psi_N|^2)\psi_N,
		\end{aligned}
		& \mathbf{x}\in\Omega,\ t>0,\\
		\psi_N(\mathbf{x},0)=\psi_0(\mathbf{x}),
		& \mathbf{x}\in\overline{\Omega},
		\end{cases}
	\end{equation}
	where $\psi_N:=\psi_N(\mathbf{x},t)$. Formally, the NLSEwPP \cref{eq::nlsem} preserves the energy
	    \begin{equation}\label{eq:EN}
	    \begin{aligned}
	        E_N(\psi_N(\cdot, t))={}&\int_{\Omega} \left[|\nabla \psi_N|^2
	        +(P_N V)(\vx)|\psi_N|^2 +\frac{\beta}{2} |\psi_N|^4\right] \rmd \vx.
	    \end{aligned}
	    \end{equation}    
    We have the following estimate for the difference of the energy $E$ and $E_N$. 
    
    \begin{proposition}[Energy error]\label{prop:energy_error}
        Assuming that $V \in H^\alpha_\text{\rm per}(\Omega)$ for some $\alpha \geq 0$ and $\psi_0 \in H_\text{\rm per}^m(\Omega)$ for some $m \geq 2$, we have
        \begin{equation*}
            |E(\psi_0) - E_N(\psi_0)| \lesssim \frac{1}{N^{\alpha + m}}. 
        \end{equation*}
    \end{proposition}

    \begin{proof}
        Recalling \cref{eq:energy,eq:EN}, by the orthogonality, we have
	        \begin{align*}
	            E(\psi_0) - E_N(\psi_0)
	            &= \int_{\Omega} (V(\vx) - P_N V(\vx)) |\psi_0(\vx)|^2 \rmd \vx \\
	            &= \langle (I-P_N)V,|\psi_0|^2\rangle
	            = \langle (I-P_N)V,(I-P_N)|\psi_0|^2\rangle,
	        \end{align*}
        where $\langle \cdot, \cdot \rangle$ is the standard $L^2$ inner-product. It follows from the Cauchy-Schwarz inequality, the standard projection error estimates, and $|\psi_0|^2 \in H^m_{\text{\rm per}}$ that
        \begin{equation*}
            |E(\psi_0) - E_N(\psi_0)| \lesssim \| (I-P_N) V \|_{L^2} \| (I-P_N)|\psi_0|^2 \|_{L^2} \lesssim \frac{1}{N^{\alpha}} \frac{1}{N^m}, 
        \end{equation*}
        which completes the proof. 
    \end{proof}
    
    Hence, even for singular potential that is merely in $L^2$, the energy error is at the same order as the Fourier projection error of the initial data. Note that although the solution can only have limited regularity under low-regularity or singular potential in general, the initial data can still be smooth. 

\subsection{Well-posedness and convergence}
We next turn to the well-posedness of the NLSEwPP \cref{eq::nlsem} and its convergence to the NLSE \cref{eq::nlse}. 

To prove the $H^2$ well-posedness of the NLSEwPP  \cref{eq::nlsem}, we make the following assumptions on the potential $V$ and the solution $\psi$ to the NLSE \cref{eq::nlse}. 

\noindent\textbf{Assumptions:} 
% (i)Assume that $\sigma= 1$ in the power-type nonlinearity \cw{Not an assumption. (i) can be removed}
%     \[
%         f(|\psi|^2)\psi = \beta |\psi|^{2}\psi .
%     \]
%     Then the nonlinear mapping $z\mapsto \beta |z|^{2}z$ is locally Lipschitz continuous. More precisely, for any $M_0>0$, there exists a constant $C_{\mathrm{Lip}}(M_0)>0$ such that
%     \begin{equation}
%     \label{eq::nona}
%         \left| f(|z_1|^{2})z_1
%         - f(|z_2|)^{2}z_2 \right|
%         \le C_{\mathrm{Lip}}(M_0)|z_1-z_2|,
%         \quad z_j\in\mathbb{C},\ |z_j|\le M_0,\ j=1,2.
%     \end{equation}
(i) Assume that $V\in L^\infty(\Omega)$ satisfies 
\begin{equation}\label{eq:PNV_bound}
	M_V:=\sup_{N \geq 1} \| P_N V \|_{L^\infty} < \infty.
\end{equation}
(ii) Let $0<T<T_\text{max}$ with $T_\text{max}$ being the maximal existence time of the solution to the NLSE \eqref{eq::nlse}. Assume that $\psi \in C([0, T]; H^2_\text{per}(\Omega)) \cap C^1([0, T]; L^2(\Omega))$ is the unique solution to \eqref{eq::nlse} with
\begin{equation}
	M_\psi:=\sup_{0 \leq t \leq T} \| \psi(t) \|_{H^2} < \infty.  
\end{equation} 

As mentioned before, \cref{eq::nlsem} makes sense for all $V \in L^2(\Omega)$. However, we need a stronger assumption $V \in L^\infty(\Omega)$ and $P_N V \in L^\infty(\Omega)$ uniformly in $N$ to prove the (conditional) $H^2$ convergence of $\psi_N$ to $\psi$ in \cref{lemma::dt_psiN_H2}, which will allow us to iteratively extend the existence time of $\psi_N$ to the same as $\psi$. Nevertheless, the proposed PP-FP method introduced in \cref{sec::VP} works well even for $L^2$-potentials, as demonstrated by the numerical results in \cref{sec:4}.    

\begin{remark}
	Although \cref{eq:PNV_bound} does not follow from $V \in L^\infty(\Omega)$, it is satisfied by most of the physically meaningful $V \in L^\infty(\Omega)$ below:
%    The assumption \cref{eq:PNV_bound} is satisfied by several relevant classes of potentials:
		\begin{enumerate}[(i)]
			\item In one dimension (1D), the assumption \cref{eq:PNV_bound} is satisfied by $V$ of bounded variation (BV) \cite{zygmund2003}. More precisely, we have
			\begin{equation*}
				\sup_{N\geq 1}\|P_NV\|_{L^\infty}
				\lesssim (\|V\|_{L^1}+|V|_\text{TV}), 
			\end{equation*}
			where $|V|_\text{TV}$ is the total variation of $V$ on $\Omega$. Hence, bounded piecewise smooth or piecewise constant potentials with finitely many discontinuities are covered. 
			
			\item In higher dimensions, a simple sufficient condition is that $V$ is the (finite linear combination of) tensor product of one-dimensional BV functions. 
%            \begin{equation}
%                V(\vx) = \sum_{j=1}^J V_j(\vx), \quad V_j(\vx) = \prod_{k=1}^d h_{j, k}(x_k), 
%            \end{equation}
%            where each $h_{j, k}$ as a one-dimensional function is bounded and of bounded variation. Indeed, we have
%            \begin{equation}
%                (P_N V_j)(\vx) = \prod_{k=1}^d (P_N^{(1)}h_{j, k})(x_k), 
%            \end{equation}
%            where $P_N^{(1)}$ is defined as $P_N$ with $d=1$. 
            Typical examples include finite linear combinations of characteristic functions of rectangles in two dimensions (2D), or boxes in three dimensions (3D). 
			
			\item For continuous potentials, \cref{eq:PNV_bound} is satisfied if $V$ satisfies
			\begin{equation*}
				\sum_{\mathbf k \in\mathbb Z^d}|\widehat V_{\mathbf k}|<\infty. 
			\end{equation*}
			In particular, this covers $V\in H_\text{per}^\alpha$ with $\alpha>d/2$.
		\end{enumerate}
		In fact, we can also remove the assumption \cref{eq:PNV_bound} if defining $P_N$ as a smooth cut-off in the Fourier space. Since this modification does not make a significant difference in practice, we will not discuss it in detail. 
\end{remark}

% We have \cref{eq:PNV_bound} is satisfied by the following potential: 
% \begin{itemize}[(i)]
% 	\item power-law trapping potentials on a bounded domain: \[ V(\vx)=|\vx|^\gamma,\quad \gamma>0,\quad \vx\in\Omega ; \]
%     \item square-well, barrier, or step potentials: \[ V(\vx)=\sum_{j=1}^J V_j\chi_{\Omega_j}(\vx), \quad J<\infty,\quad V_j\in\mathbb R ; \]
% \end{itemize}

Then we have the following result concerning the well-posedness of the NLSEwPP \eqref{eq::nlsem} and its convergence.
\begin{theorem}\label{thm::error_regularization}
    Under the assumptions (i)-(ii), there exists $N_0 \geq 1$ sufficiently large such that for all $N \geq N_0$, \eqref{eq::nlsem} admits a unique solution $\psi_N$ on the time interval $[0, T]$ satisfying $\psi_N \in C([0, T]; H^2_{\rm per}(\Omega)) \cap C^1([0, T]; L^2(\Omega))$, and 
    \begin{equation}\label{eq:uniformH2}
    	\sup_{N \geq N_0}\| \psi_N \|_{L^\infty([0, T]; H^2)} \leq M_\psi + 1.
    \end{equation}
    Moreover, we have
    \begin{equation}\label{eq:convergenceL2H1}
        \| \psi - \psi_N \|_{L^\infty([0, T];L^2)} \lesssim \frac{1}{N^2}, \quad \| \psi - \psi_N \|_{L^\infty([0, T];H^1)} \lesssim \frac{1}{N}.
    \end{equation} 
\end{theorem}

As a corollary of \cref{thm::error_regularization}, if the potential and the initial data have better regularity, we can get better regularity for $\psi$ and $\psi_N$, and faster convergence rates. 
\begin{corollary}[Higher regularity]\label{cor:high-regularity}
    Under the same assumptions as \cref{thm::error_regularization}, if, in addition, $V \in H_\text{\rm per}^\alpha(\Omega)$ with $\alpha >d/2$ and $\psi_0 \in H_\text{\rm per}^{2+\alpha}(\Omega)$, we have
    \begin{equation}\label{eq:uniform_bound_high}
        \| \psi \|_{L^\infty([0, T]; H^{2+\alpha})} < \infty, \quad \sup_{N \geq N_0} \| \psi_N \|_{L^\infty([0, T]; H^{2+\alpha})} < \infty, 
    \end{equation}
    and
    \begin{equation}\label{eq:convergenceL2H1high}
        \| \psi - \psi_N \|_{L^\infty([0, T];L^2)} \lesssim \frac{1}{N^{2+\alpha}}, \quad \| \psi - \psi_N \|_{L^\infty([0, T];H^1)} \lesssim \frac{1}{N^{1+\alpha}}.
    \end{equation}
\end{corollary}

% We first state the main convergence result. It shows that, on any compact time interval contained in the maximal existence interval of the original NLSE, the solution $\psi_N$ of the equation with potential projection converges to the solution $\psi$ of the original NLSE. Moreover, the convergence rates are of order $O(N^{-2})$ in the $L^2$ norm and $O(N^{-1})$ in the $H^1$ norm.

% Let $\psi_N$ and $\psi$ be the solutions to \eqref{eq::nlsem} and \eqref{eq::nlse}, respectively, with the same initial data
% \[
%     \psi_N(0)=\psi(0)=\psi_0 .
% \]
% Then we have Theorem \ref{thm::convmr}.

% \begin{theorem}
% \label{thm::convmr}
% Let $t_{\max}$ be the maximal existence time of the solution $\psi$ to the original NLSE \eqref{eq::nlse}. Assume that, for any $T\in(0,t_{\max})$,
% \[
%     \psi \in C([0,T],H^2_{\mathrm{per}}(\Omega))
%     \cap C^1([0,T],L^2(\Omega)).
% \]
% Then there exists $N_0>0$ such that, for all $N\ge N_0$, the equation with potential projection \eqref{eq::nlsem} admits a solution
% \[
%     \psi_N\in C([0,T],H^2_{\mathrm{per}}(\Omega))
%     \cap C^1([0,T],L^2(\Omega)).
% \]
% Moreover, for all $0\le t\le T$, the following estimates hold:
% \begin{equation}
% \label{eq::ltre} 
%     \|\psi(\cdot,t)-\psi_N(\cdot,t)\|_{L^2}
%     \lesssim N^{-2}, 
%     \quad
%     \|\psi(\cdot,t)-\psi_N(\cdot,t)\|_{H^1}
%     \lesssim N^{-1}.
% \end{equation}
% \end{theorem}

\begin{remark}
    The (optimal) convergence rates in \cref{eq:convergenceL2H1,eq:convergenceL2H1high} are essential in the design of the TS-PP-FP methods in \cref{sec::VP} since the projection parameter $N$ also determines the degree of freedom (DOF) needed in the FP discretization, as is shown in \cref{sec:4}.
\end{remark}

To prove Theorem \ref{thm::error_regularization}, we first need the following local well-posedness of the $H^2$-solution to the NLSE. 

\begin{proposition}[$H^2$ local theory]\label{prop::local_existence}
    Let $R_1, R_2 > 0$. There exist $T_0, M_0>0$ depending on $R_1, R_2$ such that whenever $\| W \|_{L^\infty} \leq R_1$ and $\| u_0 \|_{H^2} \leq R_2$, the NLSE 
    \[
    i\partial_tu(t)=-\Delta u(t) + Wu(t) +f(|u(t)|^2)u(t),\quad u(0)=u_0,
    \]
    has a unique solution $u\in C([0,T_0];H^2_{\mathrm{per}}(\Omega)) \cap C^1([0,T_0];L^2(\Omega))$
    with
    \begin{equation}\label{eq:uniform_bound}
    	\sup_{0 \leq t \leq T_0} (\|\partial_t u(t) \|_{L^2} + \| u(t) \|_{H^2} ) \leq M_0.  
    \end{equation}
\end{proposition}

\begin{proof}
	The result can be proved by the same argument used in the proof of Theorem 4.8.1 in Cazenave~\cite{cazenave2003semilinear} for the local $H^2$ well-posedness on $\R^d$. Here, under the stronger assumption of $L^\infty$-potential, we can avoid the use of the Strichartz estimates and thus extend the proof to the periodic setting directly. We omit the details for brevity. 
%	When $\Omega = \mathbb{R}^d$, with the use of Strichartz estimates, the $H^2$ local well-posedness is established for $V \in L^2$ (see, for instance). The key technique to overcome the low regularity of the potential is to use Kato's idea: to estimate the $H^2$-norm, instead of differentiating the equation twice in space, one can differentiate the equation once in time, and induce the control of the $H^2$-norm from the equation. The proof proceeds by performing a fixed-point argument on a properly chosen auxiliary space ((4.8.9) in \cite{cazenave2003semilinear}), and the estimate in \cref{eq:uniform_bound} follows directly from the construction of that space. 
%	
%	With \cref{eq:uniform_bound}, \cref{eq:H2_regularity} can be obtained from \cref{eq::nlse} directly: from \cref{eq::nlse}, taking the $L^2$-norm, using the Sobolev embedding $H^\frac{7}{4} \hookrightarrow L^\infty$, we get
%	\begin{equation*}
%		\| \Delta \psi(t) \|_{L^2} \leq \| \partial_t \psi(t) \|_{L^2} + \| V \psi(t) \|_{L^2} + \| f(|\psi(t)|^2)\psi(t) \|_{L^2} \lesssim M_0 + \| V \|_{L^\infty} M_0 + C (M_0). 
%	\end{equation*}
%	The proof is thus completed by taking the supremum over $t \in [0, T_0]$. 
\end{proof}

With \cref{prop::local_existence}, we can apply it to \cref{eq::nlsem} to get local solutions $\psi_N$ on the time interval $[0, T_0]$ with the estimates \cref{eq:uniform_bound}. Note that both $T_0$ and $M_0$ are uniform in $N$ under the assumption \cref{eq:PNV_bound}. To iteratively apply \cref{prop::local_existence} to extend the local solutions to the whole time interval $[0, T]$, we need to control the $H^2$-norm of $\psi_N$, which is a consequence of the convergence of $\psi_N$ to $\psi$ in the $H^2$-norm. 

To establish the convergence of $\psi_N$ to $\psi$ in the $H^2$-norm, we need the following. 

\begin{lemma}\label{lem:multiplier}
	%Suppose $a_N\in L^\infty(\Om)$ satisfies
	%\[
	%    \sup_N\norm[L^\infty]{a_N}<\infty,
	%    \quad
	%    \norm[L^2]{a_N}\longrightarrow0.
	%\]
	%Then $a_Ng\to0$ in $L^2$ for every $g\in L^2$.  Furthermore, if
	%$\mathcal K\subset L^2$ is compact, then
	%\begin{equation}
	%    \sup_{g\in\mathcal K}\norm[L^2]{a_Ng}\longrightarrow0.
	%\label{eq:uniform-compact-mult}
	%\end{equation}
	Under assumptions (i) and (ii), we have
	\begin{align}
		&\lim_{N \rightarrow \infty} \delta_1(N) = 0, \quad  \delta_1(N):=\sup_{0 \leq t \leq T} \| (V - P_N V) \psi(t) \|_{L^2}, \label{lem:delta1} \\
		&\lim_{N \rightarrow \infty} \delta_2(N) = 0, \quad \delta_2(N):=\sup_{0 \leq t \leq T} \| (V - P_N V) \partial_t \psi(t) \|_{L^2}. \label{lem:delta2}
	\end{align}
\end{lemma}

\begin{proof}
	By the assumption (i), we have
	\begin{equation*}
		\lim_{N \rightarrow \infty} \| V - P_N V \|_{L^2} = 0, \quad \sup_{N \geq 1} \| V - P_N V \|_{L^\infty}  \leq \| V \|_{L^\infty} + M_V < \infty. 
	\end{equation*}
Recalling $\psi \in C([0, T]; H_\text{per}^2(\Omega))$, by the Sobolev 
embedding $H^2 \hookrightarrow L^\infty$ we have
\begin{equation}
\begin{split}
\| (V - P_N V) \psi(t) \|_{L^2} &\leq \| V - P_N V \|_{L^2} \| \psi(t) \|_{L^\infty}\\
 &\lesssim \| V - P_N V \|_{L^2} \| \psi \|_{L^\infty([0, T]; H^2)} \rightarrow 0, 
 \end{split}
	\end{equation}
	which proves \eqref{lem:delta1}. To prove \eqref{lem:delta2}, we note that, for any $\phi \in L^2(\Omega)$,
	\begin{align}\label{Vdiffop}
		\| (V - P_N V) \phi \|_{L^2} 
		&\leq \| (V - P_N V) \phi 1_{|\phi| \leq K} \|_{L^2} + \| (V - P_N V) \phi (1-1_{|\phi| \leq K}) \|_{L^2} \notag \\
		&\leq  \| V - P_N V \|_{L^2} K + (\| V \|_{L^\infty} + M_V)\| \phi (1-1_{|\phi| \leq K}) \|_{L^2}, 
	\end{align}
	where $1_{|\phi| \leq K}$ is the indicator function of $\{\vx \in \Omega \ | \  |\phi(\vx)| \leq K\}$. Since $\| \phi (1-1_{|\phi| \leq K}) \|_{L^2} \rightarrow 0$ as $K \rightarrow \infty$, it follows from \eqref{Vdiffop} that
	\begin{equation}
		\lim_{N \rightarrow \infty} \| (V - P_N V) \phi \|_{L^2} = 0, \quad \phi \in L^2(\Omega). 
	\end{equation}
	This proves the convergence for any fixed $\phi \in L^2(\Omega)$. Then \eqref{lem:delta2} follows from the continuity of $\partial_t \psi(t)$, which implies that $\{\partial_t \psi(t) \ | \ t \in [0, T]\}$ is compact in $L^2(\Omega)$, and the uniform bound $ \sup_{N \geq 1} \| V - P_N V \|_{L^\infty} <\infty$. 
	%Let $C_a=\sup_N\norm[L^\infty]{a_N}$.  Given $g\in L^2$, choose
	%$g_m\in L^2\cap L^\infty$ with $g_m\to g$ in $L^2$.  Then
	%\[
	%    \norm[L^2]{a_Ng}
	%    \le C_a\norm[L^2]{g-g_m}
	%       +\norm[L^2]{a_N}\norm[L^\infty]{g_m}.
	%\]
	%First choose $m$ large and then let $N\to\infty$.  This proves strong
	%convergence on each fixed $g\in L^2$.
	%
	%The multiplication operators $M_{a_N}:g\mapsto a_Ng$ are uniformly
	%bounded on $L^2$, with $\norm{M_{a_N}}\le C_a$.  Given a compact set
	%$\mathcal K$ and $\varepsilon>0$, choose a finite
	%$\varepsilon/(2C_a)$-net $g_1,\ldots,g_J$ for $\mathcal K$.  Strong
	%convergence at these finitely many points gives
	%$\max_j\norm[L^2]{a_Ng_j}<\varepsilon/2$ for all sufficiently large
	%$N$.  Approximating an arbitrary $g\in\mathcal K$ by one of the $g_j$
	%then proves \eqref{eq:uniform-compact-mult}.
\end{proof}

\begin{proposition}[$H^2$-convergence]\label{lemma::dt_psiN_H2}
	Assume there exist $0 < T' \leq T$ and $N_0' \geq 1$ such that, for all $N \geq N_0'$, $\psi_N \in C([0, T']; H_\mathrm{per}^2(\Omega)) \cap C^1([0, T']; L^2(\Omega))$ with
	\begin{equation*}
		\sup_{N \geq N_0'} \| \psi_{N} \|_{L^\infty([0, T']; H^2)} < \infty.  
	\end{equation*}
	Then we have
	\begin{equation*}
		\lim_{N \rightarrow \infty} \| \psi - \psi_{N} \|_{L^\infty([0, T']; H^2)} = 0. 
	\end{equation*}
\end{proposition}

\begin{proof}
	The assumption $\sup_{N \geq N_0'} \| \psi_N \|_{L^\infty([0, T']; H^2)} < \infty$ also implies that
	\begin{equation}
		\sup_{N \geq N_0'} \| \partial_t \psi_N \|_{L^\infty([0, T']; L^2)} < \infty,  
	\end{equation}
	which follows from \cref{eq::nlsem} through
	\begin{equation}
		\| \partial_t \psi_N(t) \|_{L^2} \leq  \| \psi_N(t) \|_{H^2} + \| P_N V \|_{L^\infty} \| \psi_N(t) \|_{L^2} + C(\| \psi_N(t) \|_{L^\infty}). 
	\end{equation}
	In the following, we define
	\begin{equation}\label{eq:R_def}
		R:=\sup_{N \geq N_0'} \| \partial_t \psi_N \|_{L^\infty([0, T']; L^2)} + \sup_{N \geq N_0'} \| \psi_N \|_{L^\infty([0, T']; H^2)} < \infty. 
	\end{equation}
	We assume $N \geq N_0'$ throughout the rest of the proof. 
	
	Let $e_N(t):=\psi(t) - \psi_N(t)$ for $0 \leq t \leq T'$. Subtracting \cref{eq::nlsem} from \cref{eq::nlse}, we obtain
	\begin{equation}\label{eq:eN}
		i \partial_t e_N(t) = -\Delta e_N(t) + (V - P_N V)\psi(t) + P_N V e_N(t) + F(\psi(t)) - F(\psi_N(t)), 
	\end{equation}
	where we denote $F(\phi) = f(|\phi|^2)\phi$. Hence, by the Duhamel's formula for \cref{eq:eN}, we have, noting that $e_N(0) = 0$, 
	\begin{equation}\label{eq:eNduhamel}
		e_N(t) = -i \int_0^t e^{i(t-s)\Delta} ((V - P_N V) \psi(s) + P_N V e_N(s) + F(\psi(s)) - F(\psi_N(s))) \rmd s. 
	\end{equation}
	Taking the $L^2$-norm on both sides, we get, recalling \cref{lem:delta1} and using $H^2 \hookrightarrow L^\infty$, 
	\begin{align}
		\| e_N(t) \|_{L^2} 
		&\leq t \sup_{0 \leq s \leq T} \| (V - P_N V) \psi(s) \|_{L^2} + \| P_N V \|_{L^\infty} \int_0^t \| e_N(s) \|_{L^2} \rmd s \notag \\
		&\quad+ \int_0^t \| F(\psi(s)) - F(\psi_N(s)) \|_{L^2} \rmd s \notag \\
		&\leq T'\delta_1(N) + C(M_V, M_\psi, R) \int_0^t \| e_N(s) \|_{L^2} \rmd s. 
	\end{align}
	By the Gronwall's inequality and \cref{lem:delta1}, we have
	\begin{equation}\label{eq:convergenceL2}
		\lim_{N \rightarrow \infty} \| e_N \|_{L^\infty([0, T']; L^2)}  = 0. 
	\end{equation}
	Note that $\sup_{N \geq N_0'} \| e_N \|_{L^\infty([0, T']; H^2)} < \infty$, we have, by the Sobolev embedding and the interpolation inequality, 
	\begin{equation}
		\| e_N(t) \|_{L^\infty} \lesssim \| e_N(t) \|_{H^\frac{7}{4}} \lesssim \| e_N(t) \|_{L^2}^\frac{1}{8}  \| e_N(t) \|_{H^2}^\frac{7}{8}. 
	\end{equation}
	As a result, 
	\begin{equation}\label{eq:convergenceLinfty}
		\lim_{N \rightarrow \infty} \| e_N \|_{L^\infty([0, T']; L^\infty)}  = 0. 
	\end{equation}
	Then, to show the convergence of $e_N(t)$ in the $H^2$-norm, it suffices to prove the convergence of $\| \partial_t e_N(t) \|_{L^2}$ and use \cref{eq:eN} to recover the $H^2$ convergence through
	\begin{equation}\label{eq:convergenceDelta}
		\| \Delta e_N(t) \|_{L^2} \leq \| \partial_t e_N(t) \|_{L^2} + \delta_1(N) + C(M_V, M_\psi, R) \| e_N(t) \|_{L^2}. 
	\end{equation}
	
	To estimate $\partial_t e_N$, we rewrite \cref{eq:eNduhamel} as
	\begin{equation*}
		e_N(t) = -i \int_0^t e^{is\Delta} ((V - P_N V) \psi(t-s) + P_N V e_N(t-s) + F(\psi(t-s)) - F(\psi_N(t-s))) \rmd s. 
	\end{equation*}
	Taking the time derivative and changing the variable back, we obtain
	\begin{align}\label{eq:dt_duhamel}
		\partial_t e_N(t) 
		&= -ie^{it\Delta} (V - P_N V) \psi_0 -i \int_0^t e^{i(t-s)\Delta} \Bigl[(V - P_N V) \partial_s \psi(s)
		+ P_N V \partial_s e_N(s) \notag\\
		&\hspace{13em}{}+ \partial_s \left(F(\psi(s)) - F(\psi_N(s))\right)\Bigr] \rmd s,
	\end{align}
	where the differentiation under the integral is justified by noting that for any $u \in C([0, T']; H^2(\Omega)) \cap C^1([0, T']; L^2(\Omega))$, $F(u) \in C^1([0, T']; L^2(\Omega))$ with
	\begin{equation}
		\partial_t F(u(t)) = \beta (u(t)^2 \overline{\partial_t u(t)} + 2|u(t)|^2 \partial_t u(t)). 
	\end{equation}
	Taking the $L^2$-norm on both sides and recalling \cref{lem:delta2}, using
	\begin{equation*}
		\| \partial_s (F(\psi(s)) - F(\psi_N(s))) \|_{L^2} \leq C(M_\psi, R) (\| \partial_s e_N(s) \|_{L^2} + \| e_N(s) \|_{L^\infty}), 
	\end{equation*} 
	we have
	\begin{equation}
	\begin{aligned}
		\| \partial_t e_N(t) \|_{L^2}
		\leq{}& \delta_1(N)+t\delta_2(N)
		+ tC(M_\psi, R)\sup_{0 \leq s \leq T'} \| e_N(s) \|_{L^\infty}\\
		&+ C(M_V, M_\psi, R) \int_0^t \| \partial_s e_N(s) \|_{L^2} \rmd s.
	\end{aligned}
	\end{equation}
	By the Gronwall's inequality, using \cref{lem:delta2,eq:convergenceLinfty}, we have
	\begin{equation}\label{eq:convergencedt}
		\lim_{N \rightarrow \infty} \| \partial_t e_N \|_{L^\infty([0, T']; L^2)}  = 0. 
	\end{equation}
	By \cref{eq:convergenceDelta}, together with \cref{eq:convergencedt,eq:convergenceL2}, we conclude that
	\begin{equation}
		\lim_{N \rightarrow \infty} \| e_N \|_{L^\infty([0, T']; H^2)}  = 0. 
	\end{equation}
	The proof is thus completed. 
\end{proof}

Though we have established the convergence in the $L^2$-norm in the proof of \cref{lemma::dt_psiN_H2}, we have the following result that refines the convergence rate.

Recall the following property from \cite{bao2024efp,bao2024ewinls}: For every $\phi\in L^2(\Omega)$ and $W \in L^\infty(\Omega)$,
\begin{equation}\label{lem:frequency-separation}
    P_{N/2}((W - P_N W) P_{N/2} \phi)=0.
\end{equation}

\begin{proposition}[$L^2$-convergence]\label{lemma::dt_psiN_L2}
%	Assume there exists $N_0 \geq 1$ such that, for all $N \geq N_0$, $\psi_N \in C([0, T]; H_\mathrm{per}^2(\Omega)) \cap C^1([0, T]; L^2(\Omega))$ with
%	\begin{equation*}
%		\sup_{N \geq N_0} \| \psi_{N} \|_{L^\infty([0, T]; H^2)} < \infty.  
%	\end{equation*}
	Under the same assumption as \cref{lemma::dt_psiN_H2}, we have, for all $N \geq N_0'$, 
	\begin{equation*}
		\| \psi - \psi_{N} \|_{L^\infty([0, T']; L^2)} \lesssim \frac{1}{N^2},  
	\end{equation*}
	where $N_0'$ and $T'$ are defined in the assumption of \cref{lemma::dt_psiN_H2}. 
\end{proposition}

\begin{proof}
	We decompose $\psi(t) - \psi_N(t)$ as
	\begin{equation}\label{eq:decomp}
		\psi(t) - \psi_N(t) = P_{N/2} (\psi(t) - \psi_N(t)) + (I-P_{N/2})\psi(t) - (I-P_{N/2}) \psi_N(t).
	\end{equation}
		
	Define $w_N(t): = P_{N/2} (\psi(t) - \psi_N(t))$ for $0 \leq t \leq T'$, and, again, denote $F(\phi) = f(|\phi|^2)\phi$. Applying $P_\frac{N}{2}$ to both sides of \cref{eq::nlse,eq::nlsem}, we obtain
	\begin{align}
		&i \partial_t P_\frac{N}{2} \psi(t) = -\Delta P_\frac{N}{2} \psi(t) + P_\frac{N}{2} (V \psi(t)) + P_\frac{N}{2} F(\psi(t)), \\
		&i \partial_t P_\frac{N}{2} \psi_N(t) = -\Delta P_\frac{N}{2} \psi_N(t) + P_\frac{N}{2} (P_N V \psi_N(t)) + P_\frac{N}{2} F(\psi_N(t)). 
	\end{align}
	Subtracting gives
	\begin{equation}\label{eq:wN}
		i \partial_t w_N(t) = -\Delta w_N(t) + P_\frac{N}{2}  \mathcal{T}_1 + P_\frac{N}{2}  \mathcal{T}_2,
	\end{equation}
	where
	\begin{equation}
		\mathcal{T}_1 = \mathcal{T}_1(t) = V \psi(t) - P_N V \psi_N(t), \quad \mathcal{T}_2 = \mathcal{T}_2(t) =F(\psi(t)) - F(\psi_N(t)). 
	\end{equation}
	Multiplying $\overline{w_N(t)}$ on both sides of \cref{eq:wN}, taking the imaginary part, and integrating over $\Omega$, we obtain
	\begin{equation}\label{eq:w_N_L2}
		\frac{\rmd}{\rmd t} \| w_N(t) \|_{L^2} \leq \| P_\frac{N}{2} \mathcal{T}_1(t) \|_{L^2} + \| \mathcal{T}_2(t)\|_{L^2}. 
	\end{equation}
	It remains to estimate the two terms on the RHS. 
	
	We have
	\begin{align}\label{eq:T1_decomp}
		\mathcal{T}_1 
		&=  (V - P_N V) P_\frac{N}{2} \psi(t) + P_N V w_N(t) + (V - P_N V) (I-P_\frac{N}{2}) \psi(t) \notag\\
		&\quad + P_N V (I-P_\frac{N}{2})(\psi(t) - \psi_N(t))=: \mathcal{T}_{1, 1} + \mathcal{T}_{1, 2} + \mathcal{T}_{1, 3} + \mathcal{T}_{1, 4}, 
	\end{align}
	where we note that $P_{N/2} \mathcal{T}_{1, 1} = 0$ by \cref{lem:frequency-separation}. Similarly, 
	\begin{align}\label{eq:T2_decomp}
		\mathcal{T}_2
		&=  (F(P_\frac{N}{2}\psi(t)) - F(P_\frac{N}{2}\psi_N(t))) + (F(\psi(t)) - F(P_\frac{N}{2}\psi(t))) \notag\\
		&\quad - (F(\psi_N(t)) - F(P_\frac{N}{2}\psi_N(t))) =: \mathcal{T}_{2, 1} + \mathcal{T}_{2, 2} + \mathcal{T}_{2, 3}. 
	\end{align}
%	By Duhamel's formula for \cref{eq:e_N}, we get
%	\begin{align}\label{eq:e_N_duhamel}
%		w_N(t) = e^{i t \Delta} e_N(0) -i \int_0^t e^{i(t-s)\Delta} P_\frac{N}{2} (\mathcal{T}_1(s) + \mathcal{T}_2(s)) \rmd s.
%	\end{align}
	Recalling \cref{eq:T1_decomp} and $P_{N/2} \mathcal{T}_{1, 1} = 0$, we obtain, using the projection estimates and the $H^2$-norm bounds for $\psi$ and $\psi_N$, 
	\begin{equation}\label{eq:T1_L2_est}
		\| P_\frac{N}{2} \mathcal{T}_1(t) \|_{L^2} \leq \| \mathcal{T}_{1, 2}(t) \|_{L^2} + \| \mathcal{T}_{1, 3}(t) \|_{L^2} + \| \mathcal{T}_{1, 4}(t) \|_{L^2} \lesssim \| w_N(t) \|_{L^2} + N^{-2}.
	\end{equation}
    Similarly, recalling 
    % \cref{eq::nona} and 
    \cref{eq:T2_decomp}, we have
	\begin{equation}\label{eq:T2_L2_est}
		\| \mathcal{T}_2(t) \|_{L^2} \leq \| \mathcal{T}_{2, 1}(t) \|_{L^2} + \| \mathcal{T}_{2, 2}(t) \|_{L^2} + \| \mathcal{T}_{2, 3}(t) \|_{L^2} \lesssim \| w_N(t) \|_{L^2} + N^{-2},
	\end{equation}
    where the constants depend on $M_V$, $M_\psi$, and $R$ in \cref{eq:R_def}. Inserting \cref{eq:T1_L2_est,eq:T2_L2_est} into \cref{eq:w_N_L2}, we obtain
	\begin{equation}
		\frac{\rmd}{\rmd t}\|w_N(t)\|_{L^2} \lesssim N^{-2} + \|w_N(t)\|_{L^2}, 
	\end{equation}
	which implies by the Gronwall's inequality that
	\begin{equation}
		\|w_N(t)\|_{L^2} \lesssim \frac{1}{N^2}, \quad 0 \leq t \leq T'. 
	\end{equation}
	The conclusion then follows from \cref{eq:decomp}, the standard projection error estimates, and the $H^2$-norm bound of $\psi$ and $\psi_N$ immediately. 
\end{proof}

With \cref{prop::local_existence,lemma::dt_psiN_L2,lemma::dt_psiN_H2}, \cref{thm::error_regularization} can be proved by induction.

\begin{proof}[Proof of Theorem \ref{thm::error_regularization}]
For \cref{eq::nlsem}, we first verify that
\begin{equation}
	\| \psi_N(0) \|_{H^2} = \| \psi_0 \|_{H^2} \leq M_\psi+1, \quad \| P_N V \|_{L^\infty} \leq M_V, \quad N \geq 1. 
\end{equation}
Hence, we can apply \cref{prop::local_existence} to \cref{eq::nlsem} to get $T_0 = T_0(M_V, M_\psi+1)$ such that \cref{eq::nlsem} (for all $N \geq 1$) has a unique solution $\psi_N \in C([0, T_0]; H^2_\text{\rm per}(\Omega)) \cap C^1([0, T_0]; L^2(\Omega))$ with
\begin{equation*}
	\sup_{N \geq 1} \| \psi_N \|_{L^\infty([0, T_0]; H^2)} < \infty. 
\end{equation*}
This gives the solution $\psi_N$ on the time interval $[0, T_0]$. To extend the time interval to $[0, T]$, we use an induction argument.
 
Let $L:=\lceil T/T_0 \rceil$ and divide $[0, T]$ into subintervals $I_n:=[J_n, J_{n+1}]$ for $0 \leq n \leq L-1$ where $J_n := \min\{n T_0, T\}$. Note that there are only finitely many subintervals. Assume that, for some $1 \leq m \leq L-1$, there exists $N_0 \geq 1$ such that for all $N \geq N_0$,
\begin{eqnarray}\label{eq:induction}
\qquad\ \	\psi_N \in C([0, J_m]; H^2_\text{\rm per}(\Omega)) \cap C^1([0, J_m]; L^2(\Omega)), \ \sup_{N \geq N_0} \| \psi_{N} \|_{L^\infty([0, J_m]; H^2)} < \infty.  
\end{eqnarray}  
Then we can apply \cref{lemma::dt_psiN_H2} on $[0, J_m]$ to get 
\begin{equation}\label{eq:induc1}
	\lim_{N \rightarrow \infty}\| \psi - \psi_N \|_{L^\infty([0, J_m]; H^2)} = 0. 
\end{equation}
Hence, by choosing a possibly larger $N_0$, we have  
\begin{equation}\label{eq:induc2}
	\sup_{N \geq N_0} \| \psi_N(J_{m}) \|_{H^2} \leq M_\psi+1.  
\end{equation}
Applying \cref{prop::local_existence}, we can extend the solutions $\psi_N$ to $[0, J_{m+1}]$ with \cref{eq:induction} holding for $m+1$. This concludes the induction step. Applying \cref{lemma::dt_psiN_H2} for $[0, T]$ and choosing a possibly larger $N_0$ prove the first part of \cref{thm::error_regularization}. 

To prove \cref{eq:convergenceL2H1}, we directly apply \cref{lemma::dt_psiN_L2} for the time interval $[0, T]$, which gives the convergence in the $L^2$-norm. The convergence in the $H^1$-norm follows immediately by interpolation between $L^2$ and $H^2$. Then the proof for \cref{thm::error_regularization} is completed.
\end{proof}

\begin{proof}[Proof of \cref{cor:high-regularity}]
    When $d/2<\alpha \leq 2$, the proof follows the same idea as \cref{thm::error_regularization}, and can be extended to general $\alpha >2$ by a finite bootstrap argument. We briefly explain the idea of the proof for $d/2<\alpha \leq 2$ below. 
    
    We note that $\| P_N V \|_{H^\alpha} \leq \| V \|_{H^\alpha}$, and $H^\alpha$ is an algebra when $\alpha > d/2$. To obtain \cref{eq:uniform_bound_high}, we first prove the $H^\alpha$-regularity for the time derivatives $\partial_t \psi$ and $\partial_t \psi_N$ through (differentiating the Duhamel formula similar to \cref{eq:dt_duhamel}), 
    \begin{equation*}
    \begin{aligned}
        \partial_t v(t) 
        &= e^{i t \Delta} i \Delta \psi_0 - i e^{it \Delta} (W\psi_0+F(\psi_0)) \notag \\
        &\quad - i \int_0^t e^{i(t-s)\Delta} (W \partial_t v(s) + \partial_t F(v(s))) \rmd s, \quad (v, W) \in \{(\psi, V), (\psi_N, P_N V)\},
    \end{aligned} 
    \end{equation*}
    where the $H^2$-regularity of $\psi$ and $\psi_N$ are also used to control the nonlinearity. Then  one can conclude the $H^{2+\alpha}$-regularity from the equation. Since the $H^2$ solutions are already constructed on $[0,T]$, no higher-regularity local theory is needed. For \cref{eq:convergenceL2H1high}, the $(2+\alpha)$-order convergence in the $L^2$-norm is a direct consequence of the $H^{2+\alpha}$ regularity of the solutions and the standard projection error estimates, together with the frequency separation property \cref{lem:frequency-separation}. The convergence in the $H^1$-norm again follows from interpolation. We omit the details of the proof for brevity. 
\end{proof}

\subsection{Numerical validation}
In this subsection, we present numerical results to verify the proved convergence rate of $\psi_N \rightarrow \psi$ in \cref{thm::error_regularization,cor:high-regularity}. 

% \qs{$H^2$ initial data tested. The observed spatial convergence orders for the $H^1$ and $H^2$ potentials remain approximately three and four, respectively. Maybe the prefactor of the second-order term may be too small for the asymptotic second-order behavior to be visible over the tested range.}
% To validate the above theoretical findings, we present numerical experiments that examine the convergence behavior of the equations with potential projection as the projection parameter $N$ increases. In particular, we verify the rates established in \cref{thm::error_regularization} in both the $L^2$ and $H^1$ norms. We consider the cubic nonlinear Schr\"{o}dinger equation
% \begin{equation}
% \label{eq::nlsets}
%     \begin{cases}
%     i\partial_t\psi(\mathbf{x},t)
%     =-\Delta\psi(\mathbf{x},t)+V(\mathbf{x})\psi(\mathbf{x},t)
%     -|\psi(\mathbf{x},t)|^2\psi(\mathbf{x},t),
%     \quad \mathbf{x}\in\Omega,\ t>0,\\
%     \psi(\mathbf{x},0)=\psi_0(\mathbf{x})
%     =e^{-|\mathbf{x}|^2/2},
%     \quad \mathbf{x}\in\overline{\Omega}.
%     \end{cases}
% \end{equation}

% The tests are carried out in one and two spatial dimensions. In the one-dimensional examples, we take $\Omega=(-16,16)$, while in the two-dimensional examples we take $\Omega=(-8,8)^2$.
The tests are carried out in one spatial dimension with $\Omega=(-16,16)$. We choose three representative (low-regularity) potentials $V = V_j \ (j=1/2, 1, 2)$:  for $x \in \Omega$,
%\cref{eq::vs1}\cw{give the potential}. 
\begin{equation}
\label{eq::vs1}
\begin{aligned} 
V_\frac{1}{2}(x) &=
\left\{
\begin{aligned}
&-4, && x\in(-2,2), \\
&0, && \text{otherwise}, 
\end{aligned}
\right. , \quad V_{1}(x)= |x|^{0.51}, \quad V_2(x)= |x|^{1.51}\left(1 - \frac{x^2}{16^2}\right)^2.
\end{aligned}
\end{equation}
It can be checked that $V_{1/2} \in L^\infty \cap H^{\frac{1}{2}-}$, $V_1 \in H^1$, and $ V_2 \in H^2$. In addition, though not covered in \cref{thm::error_regularization}, we also test the convergence of $\psi_N \rightarrow \psi$ for a singular $L^2$-potential 
\begin{equation}\label{eq:V0}
    V_0(x)=|x|^{-0.49}. 
\end{equation}
In all the cases, we choose $\beta = 1$ and a Gaussian initial datum $\psi_0(x) = e^{-x^2/2}$. 
%which explains the observed $2.5$-order convergence in the $L^2$-norm, and the $1.5$-order convergence in the $H^1$-norm. 

\begin{figure}[htbp]
\centering
    \includegraphics[width=0.475\linewidth]{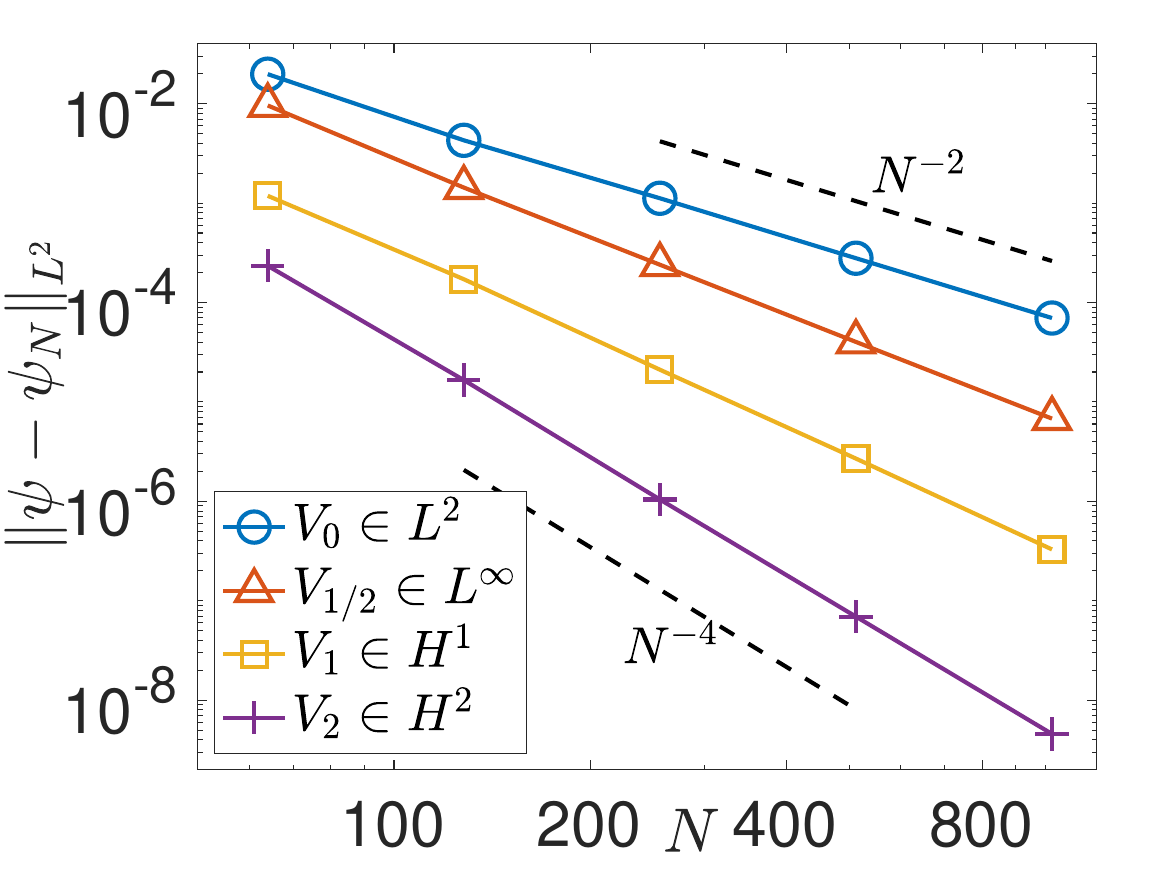}
    \includegraphics[width=0.475\linewidth]{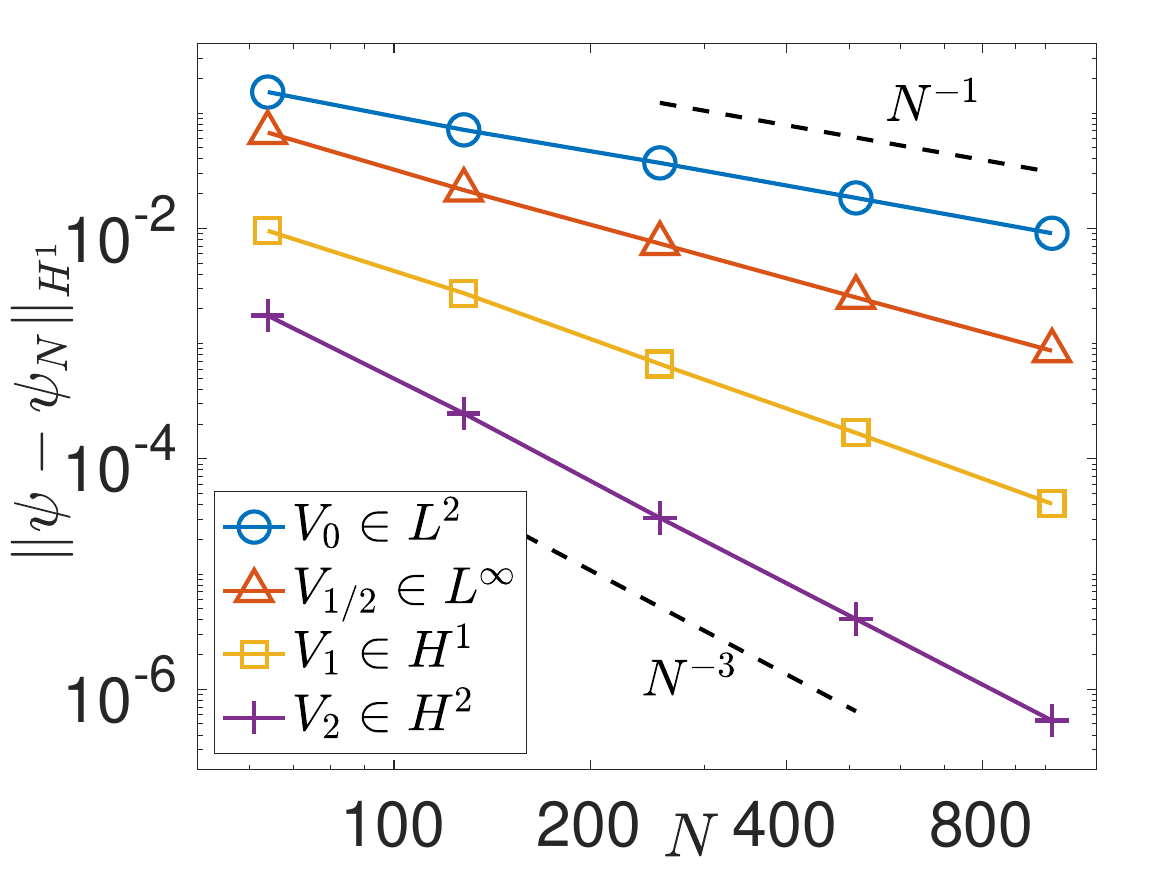}
\caption{Potential projection error $\psi - \psi_N$ in the $L^2$-norm (left) and the $H^1$-norm (right) in terms of $N$ for potentials $V_0$, $V_{1/2}$, $V_1$, and $V_2$}
\label{fig:section21}
\end{figure}

The numerical results are presented in \cref{fig:section21}, where we can clearly observe the proved convergence order in \cref{thm::error_regularization,cor:high-regularity}, indicating the optimality of our estimates. Moreover, the numerical result for $V_0$ suggests that \cref{thm::error_regularization} may be valid even for singular $L^2$-potentials. 

\section{Time-splitting PP-FP methods}\label{sec::VP}

We combine the PP-FP method with time-splitting methods to construct fully discrete numerical schemes for the NLSE \cref{eq::nlse} with low-regularity or singular potential. 

The time-splitting methods are very popular in practice due to their efficiency in the implementation and many structure-preservation properties such as the mass conservation, time symmetry, gauge invariance, preservation of dispersion relation, etc. In fact, one particular advantage of the proposed PP-FP method (compared to the FS and eFP methods \cite{bao2024efp,bao2024tsnls}) in the low regularity setting is the preservation of these structure-preservation properties inherited by the splitting methods at the time-semidiscrete level. 
% In particular, the preservation of the $L^2$-norm demonstrates the potential application of the proposed time-splitting PP-FP methods in quantum computing. 
Again, although we only consider the time-splitting methods for the time discretization, the same idea can be directly generalized to other time integrators such as exponential integrators and finite difference methods.

The time-splitting PP-FP method can be constructed by simply applying the standard TS-FP method \cite{bao2013bec} to the NLSEwPP, with a properly chosen spatial disccretization parameter $M \geq N$ depending on the projection parameter $N$. 
% The numerical results suggest that $M = N$ or more general $M \geq N$ need to be used in practice. In contrary, $M = N/2$ would totally destroy the expected convergence rate of the PP-FP method. 
Since the TS-FP method is standard, we only present the numerical scheme for the first-order Lie-Trotter splitting. The generalization to other splitting schemes is straightforward. 

\subsection{Semidiscretization by time-splitting methods}
We recall the semidiscrete time-splitting schemes in this subsection. We shall derive the numerical scheme for the NLSE \cref{eq::nlse} by directly discretizing the NLSEwPP \cref{eq::nlsem}.

The time-splitting methods are based on the decomposition of the original equation into two subproblems such that both can be solved easily. For the NLSEwPP \cref{eq::nlsem}, this decomposition reads
\begin{equation}
    \partial_t\psi=A(\psi)+B(\psi), \qquad A(\psi):=i\Delta\psi,\quad
    B(\psi):=-i\left(P_N V+f(|\psi|^2)\right)\psi.
\end{equation}
% where
% \begin{equation}\label{eq::ab}
% \end{equation}
% Here $A$ represents the linear dispersive part, while $B$ contains the external potential and the nonlinear interaction.
The flow generated by $A$ is governed by the free Schr\"{o}dinger equation
\begin{equation}
\begin{cases}
    \partial_t\psi(x,t)=i\Delta\psi(x,t),\quad x\in\Omega,\quad t>0,\\
    \psi(x,0)=\psi_0(x),\quad x\in\overline{\Omega},
\end{cases}
\end{equation}
whose exact solution can be written as
\begin{equation}
    \psi(\cdot,t)=\varphi_A^t(\psi_0):=e^{it\Delta}\psi_0(\cdot),\quad t\ge0.
\end{equation}
The flow generated by $B$ is given by
\begin{equation}
\begin{cases}
    \partial_t\psi(x,t)
    =-i\left(P_N V(x)+f(|\psi(x,t)|^2)\right)\psi(x,t),\quad x\in\Omega,\quad t>0,\\
    \psi(x,0)=\psi_0(x),\quad x\in\overline{\Omega}.
\end{cases}
\end{equation}
For this subproblem, since both $P_N V(x)$ and $f(|\psi(x,t)|^2)$ are real-valued, the modulus of the solution is conserved pointwise, namely $|\psi(x,t)|=|\psi_0(x)|$ for all $t \geq 0$. Therefore, the $B$-flow can also be integrated exactly:
\begin{equation}
\label{eq::phibt}
    \psi(x,t)
    =\varphi_B^t(\psi_0)(x)
    =\psi_0(x)
    \exp\left(-it\left(P_N V(x)+f(|\psi_0(x)|^2)\right)\right),\quad t\ge0.
\end{equation}
Based on the above decomposition, different compositions of $\varphi_A^t$ and $\varphi_B^t$ lead to splitting methods of different orders and properties. 
%In particular, we present below the first-order Lie-Trotter splitting, the second-order Strang splitting, a fourth-order partitioned Runge-Kutta splitting, and a fourth-order compact splitting schemes. 

Note that both subflows preserve the $L^2$-norm. Consequently, any splitting scheme obtained by composing these exact subflows (including all the four splitting schemes considered in this section) preserves the $L^2$-norm (mass) at the time-semidiscrete level.  

Let $\tau>0$ be the time step size and denote
\[
    t_n=n\tau,\quad n=0,1,\ldots .
\]
We use $\psi^{[n]}(\cdot)$ to denote the time-semidiscrete approximation to
$\psi(\cdot,t_n)$. The first-order Lie-Trotter splitting (LT) reads 
\begin{equation}\label{eq::semilt}
    \psi^{[n+1]}=\varphi_A^\tau \circ \varphi_B^\tau \left(\psi^{[n]}\right),\quad n\geq0, \quad \psi^{[0]} = \psi_0. 
\end{equation}
Higher-order splitting schemes can be generally constructed as 
\begin{equation}\label{eq:splitting_general}
    \psi^{[n+1]} = \prod_{j = 1}^{J} \varphi_A^{a_j \tau} \circ \varphi_B^{b_j \tau}(\psi^{[n]}), \quad a_j, b_j \in \R, \quad J \geq 1, \quad \sum_{j=1}^J a_j = \sum_{j=1}^J b_j = 1.
\end{equation}
The other splitting schemes of the form \cref{eq:splitting_general} we will consider in the paper are the second-order Strang splitting (ST) (see, e.g., \cite{bao2003numerical}), and the fourth-order partitioned Runge-Kutta splitting (S4) in \cite{S4RK}. We will also consider a fourth-order compact splitting (S4c) given in \cite{S4clinear,S4c} that falls out of the general form \cref{eq:splitting_general}. The main advantage of compact splitting schemes lies in the reduction of the intermediate stages in the high-order splitting schemes. Also, it can enter the asymptotic convergence regime with larger time steps compared to the partitioned Runge-Kutta splitting schemes. In addition, it allows more freedom in the choice of the coefficients such that no negative time steps are used in the linear flow $\varphi_A^t$, which is particularly important in the application of dissipative equations.

\begin{remark}[Structure preservation]\label{rem:structure}
All the four time-splitting methods (LT, ST, S4, and S4c) preserve the mass exactly, i.e., $\| \psi^{[n]} \|_{L^2} \equiv \| \psi^{[0]} \|_{L^2}$ for $n \geq 0$. 
%Moreover, they are all gauge invariant in the sense that when $V \rightarrow V+\alpha$ for some $\alpha \in \R$, we have $P_N V \rightarrow P_N V + \alpha$, abnd thus the numerical solutions transform as $\psi^{[n]} \rightarrow e^{- i t_n \alpha} \psi^{[n]}$, leaving the density $|\psi^{[n]}|^2$ invariant under this transformation. 
Except the first-order LT scheme \cref{eq::semilt}, all the other three splitting schemes are time symmetric in the sense that the numerical scheme remains unchanged under $n+1 \leftrightarrow n$ and $\tau \rightarrow -\tau$. 
As already mentioned, all these properties are preserved after the spatial discretization by the PP-FP method.
\end{remark}

\subsection{Full discretization}
We further discretize the LT scheme \cref{eq::semilt} in space by the Fourier pseudospectral method to obtain a full discretization scheme. The extension to the other three splitting schemes is similar and omitted. 

The full discretization scheme is obtained by simply applying the standard Fourier pseudospectral method with $M \geq N$ DOF in each dimension to the semidiscretization. We emphasize the distinct roles of the two parameters here: $N$ is the parameter in the projected potential $P_N V$, whereas $M$ is the number of Fourier modes used in the Fourier pseudospectral method. As is mentioned in the introduction and also shown in \cref{fig:section11} later, the choice $M \geq N$ is very important for the performance of the PP-FP method, and other choices such as $M=N/2$ will totally destroy the superior spatial convergence behaviour of the PP-FP method.

% \subsection{First- and second-order splitting methods}

% The Lie-Trotter splitting approximates the solution by first applying flow $\varphi_B^t(\cdot)$ over one full time step and then the linear dispersive flow $\varphi_A^t(\cdot)$ over one full time step:

% This composition is expected to obtain first-order accurate in time.

% The Strang splitting uses the symmetric composition

% and is second-order accurate in time. In contrast to the Lie-Trotter splitting, the Strang splitting is symmetric with respect to time reversal. This time symmetry is an important geometric property and generally contributes to favourable long-time numerical behaviour.

For notational simplicity, we only present the full discretization in 1D, and assume that $\Omega = (a, b)$. The generalization to higher dimensions is straightforward. Let
$M$ be a positive even integer and introduce the uniform grid
\begin{equation}
x_j=a+jh,
\quad
h=\frac{b-a}{M},
\quad
j\in{T}_M^0:=\{0,1,\ldots,M-1\}.
\end{equation}
% The corresponding Fourier index set and wave numbers are
For a periodic grid function $v = (v_0, \cdots, v_M)^T \in \mathbb C^{M+1}$ with $v_0 = v_M$, we define its discrete Fourier coefficients by
\begin{equation}
\label{eq::dft}
\widetilde v_l
=
\frac{1}{M}
\sum_{j=0}^{M-1}
v_j e^{-i\mu_l(x_j-a)},
\quad l\in {T}_M := \left\{-\frac{M}{2},\ldots,\frac{M}{2}-1\right\},
\quad
\mu_l:=\frac{2\pi l}{b-a}.
\end{equation}

% Replacing the exact subflows in \cref{eq::semilt} by their discrete counterparts gives the Lie-Trotter Fourier pseudospectral method with potential projection, denoted by LT-PP-FP:

Denote by $\psi^n_j$ the fully discrete numerical approximation to $\psi(x_j, t_n)$ for $n \geq 0$ and $j \in {T}_M^0$. By precomputation, we can get $V_{N, j} := P_N V (x_j)$ on the grid points. 

Then the fully discrete LT splitting PP-FP method (LT-PP-FP) reads
\begin{equation}
\label{eq::ltsfpm}
\begin{aligned}
    \psi_j^{(1)}
    &=
    \exp\left[-i\tau\left(V_{N, j}+f(|\psi^{n}_j|^2)\right)\right]
    \psi^{n}_j,
    \quad j\in{T}_{M}^{0},\\
    \psi^{n+1}_j
    &=
    \sum_{l\in{T}_M}
    e^{-i\tau\mu_l^2}
    \widetilde{(\psi^{(1)})}_l
    e^{i\mu_l(x_j-a)},
    \quad j\in{T}_{M}^{0},\quad n\ge0, 
\end{aligned}
\end{equation}
with $\psi^0_j = \psi_0(x_j)$ for $j \in {T}_M^0$. 

It is well-known that the FP discretization of the time-splitting methods preserves the mass conservation and the time symmetry if the corresponding semidiscrete schemes do. Hence, all the four fulldiscretization schemes enjoy the same structure-preservation properties as the corresponding semidiscretizations as discussed in Remark \ref{rem:structure}.

\section{Numerical results}\label{sec:4}
% In this section, we present numerical results for the nonlinear Schr\"{o}dinger equation (NLSE). The aim is to investigate the performance of the proposed schemes. 

% \cw{In this section, we present numerical results to examine the convergence behaviour of the PP-FP methods. We also present comaprisons with the standard FP and FS methods to demonstrate its superiority in accuracy, efficiency, and structure-preservation. Moreover, we use the PP-FP method for the spatial discretization to investigate the regularity requirement on the potential for the splitting methods.}

We present numerical results to demonstrate the superiority of the proposed TS-PP-FP methods in terms of accuracy, efficiency, and structure preservation. We first examine its convergence and compare it with the Fourier pseudospectral (FP) and Fourier spectral (FS) discretization of the time-splitting methods in \cref{sec4:convergence}. We further systematically investigate the regularity requirements on the potential of the TS-PP-FP methods in \cref{sec4:splitting}.

\subsection{Spatial convergence and comparison}\label{sec4:convergence}
In this subsection, we investigate the spatial convergence of the proposed TS-PP-FP method for potentials with different regularity, and compare it with the FP and FS discretization of the time-splitting methods. To quantify the error, we define
\begin{equation*}
    e_{L^2}(t_n):=\| \psi(t_n) - I_M \psi^n \|_{L^2}, \quad e_{H^1}(t_n):=\| \psi(t_n) - I_M \psi^n \|_{H^1}, 
\end{equation*}
where $I_M$ is the standard $M$-points Foureir interpolation operator \cite{STW2021}.

% \cw{show results for both 1D and ; also test $V \in L^2$}

We first test the spatial convergence of the PP-FP method, which is implemented by using the second-order Strang splitting in time, i.e. the ST-PP-FP method in \cref{sec::VP}. We choose $\Omega = (-16,16)$, $\beta=1$, and a Gaussian initial datum $\psi_0(x) =e^{-x^2/2}$ in the NLSE \cref{eq::nlse}. For the potentials, we use the same $V = V_j$ with $j = 0, 1/2, 1, 2$ defined in \cref{eq::vs1,eq:V0}. The final time is fixed as $T=1$. Recall that $V_0 \in L^2$, $V_{1/2}\in L^\infty \cap H^{\frac{1}{2}-}$, $V_1\in H^1$, and $V_2 \in H^2$. Since the initial data is smooth, the corresponding exact solutions admit the regularity (roughly) $H^2$, $H^{2.5}$, $H^{3}$, and $H^4$, respectively. 

% \begin{equation}
% \label{eq::nlsets_reg}
% ,\quad\mathbf{x}\in\overline{\Omega}.
% \end{equation}

% \begin{equation}
% \label{eq::vs1}
% \begin{aligned}
% V_{0}(x) = 
% \begin{cases}
% \begin{aligned}
%     & -4 & x\in(-2,2) \\
%  & 0 & \text{otherwise,} & 
% \end{aligned}
% \end{cases} , 
% \quad 
% &V_{1}(x) = |x|^{0.76},\\
% \quad
% V_2(x) = |x|^{1.51}\left(1 - \frac{x^2}{16^2}\right)^2,
% \quad
% &x \in \Omega.
% \end{aligned}
% \end{equation}

% \cap W^{1,4}(\Omega),
% \[
% V_0\in L^\infty(\Omega),\quad
% V_1\in H^1_{\mathrm{per}}(\Omega),
% % \cap W^{1,4}(\Omega),
% \quad
% V_2\in H^2_{\mathrm{per}}(\Omega).
% \]

Before testing the convergence, we examine the influence of different choices of $M$ in terms of $N$ to the PP-FP method. We choose the $L^\infty$-potential $V=V_{1/2}$, and consider three choices $N=2M$, $N=M$, and $N=M/2$ (recalling that $M$ is the DOF in space and we have $M = 32/h$). In the computation, the reference solutions are computed using the TSeFP method \cite{bao2024ewinls,bao2024efp} with $\tau=\tau_e=10^{-6}$ and $h=h_e=2^{-9}$, and the numerical solutions are computed by the ST-PP-FP method \cref{eq::ltsfpm} with $\tau = \tau_e$ and varying $h$, and with $N$ chosen accordingly. 

\begin{figure}
    \centering
    \includegraphics[width=0.475\linewidth]{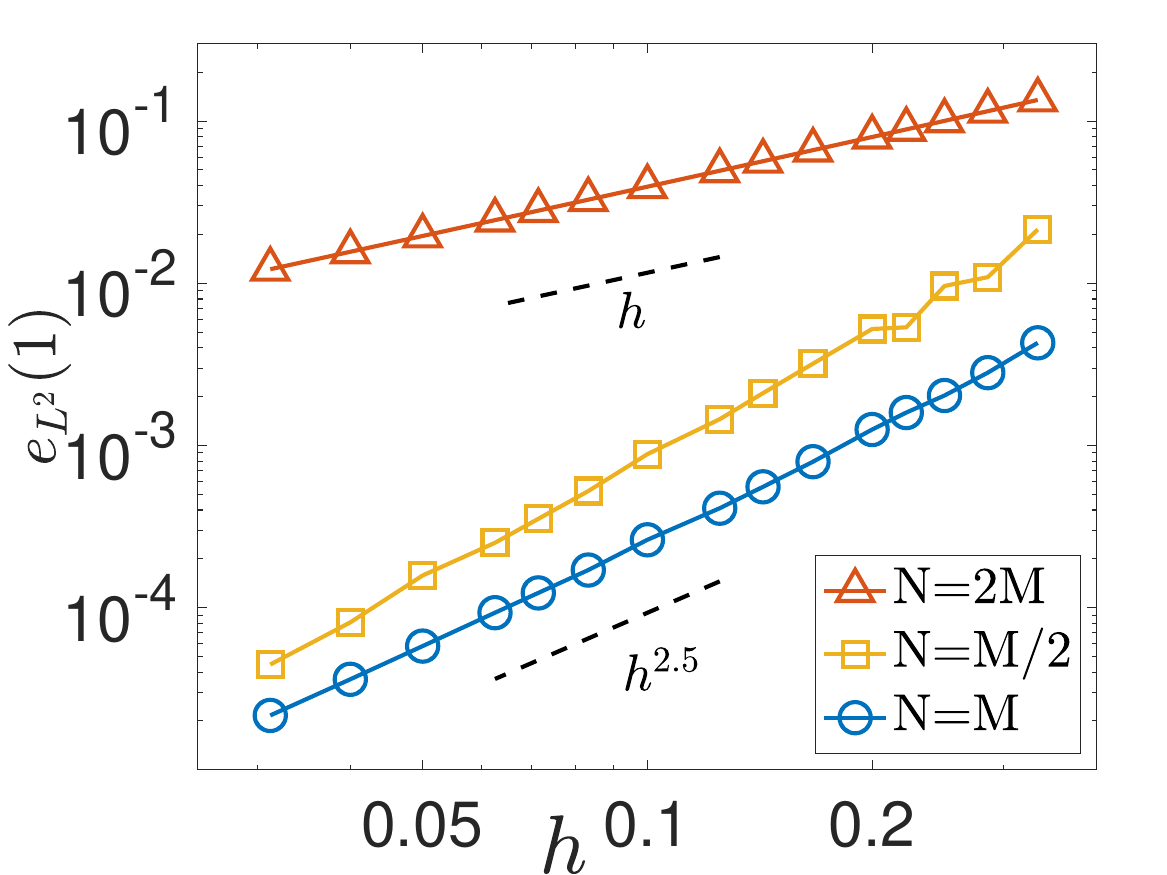}
    \includegraphics[width=0.475\linewidth]{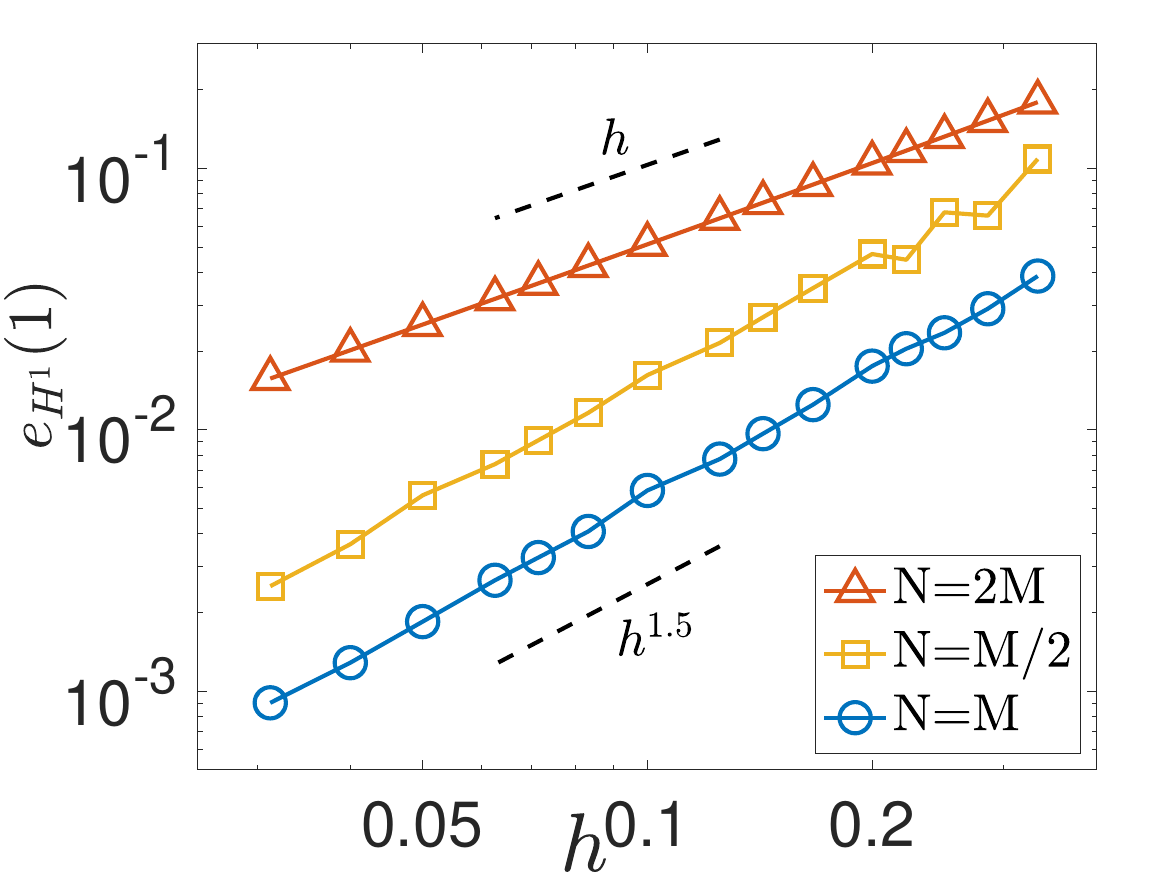}
    \caption{Errors for different choices of $N$ and $M$ in the PP-FP method: $L^2$-error (left) and $H^1$-error (right)}
    \label{fig:section11}
\end{figure}

The numerical errors are shown in \cref{fig:section11}. Among the choices considered, \(N=M\) yields the smallest errors for the PP-FP method. When \(N=2M\), the spatial grid does not adequately resolve \(P_NV\), leading to a substantial deterioration in performance: the errors are considerably larger and the observed convergence rates are significantly lower. By contrast, when \(N=M/2\), the errors lie between those obtained with \(N=M\) and \(N=2M\). They are slightly larger than those for \(N=M\), owing to the coarser potential projection, but remain much smaller than those arising from the under-resolved choice \(N=2M\). Based on this numerical observation, we see that it is necessary to choose $M \geq N$ such that the spatial discretization sufficiently resolves the projected potential $P_N V$. In practice, it is convenient to use $N=M$, which appears to give the best overall performance according to our numerical results. Therefore, in all subsequent implementations of the PP-FP method, we always take $N=M$. 

% \begin{figure}
%     \centering
%     \includegraphics[width=1\linewidth]{fig//section11.pdf}
%     \caption{Errors for different choices of $N$ and $M$ in the PP-FP method. \cw{use $L^\infty$-potential for this figure}\cw{reverse the legend order:N=2M, N=M, N=M/2}} 
%     \label{fig:section11}
% \end{figure}

\begin{figure}
    \centering
    \includegraphics[width=0.475\linewidth]{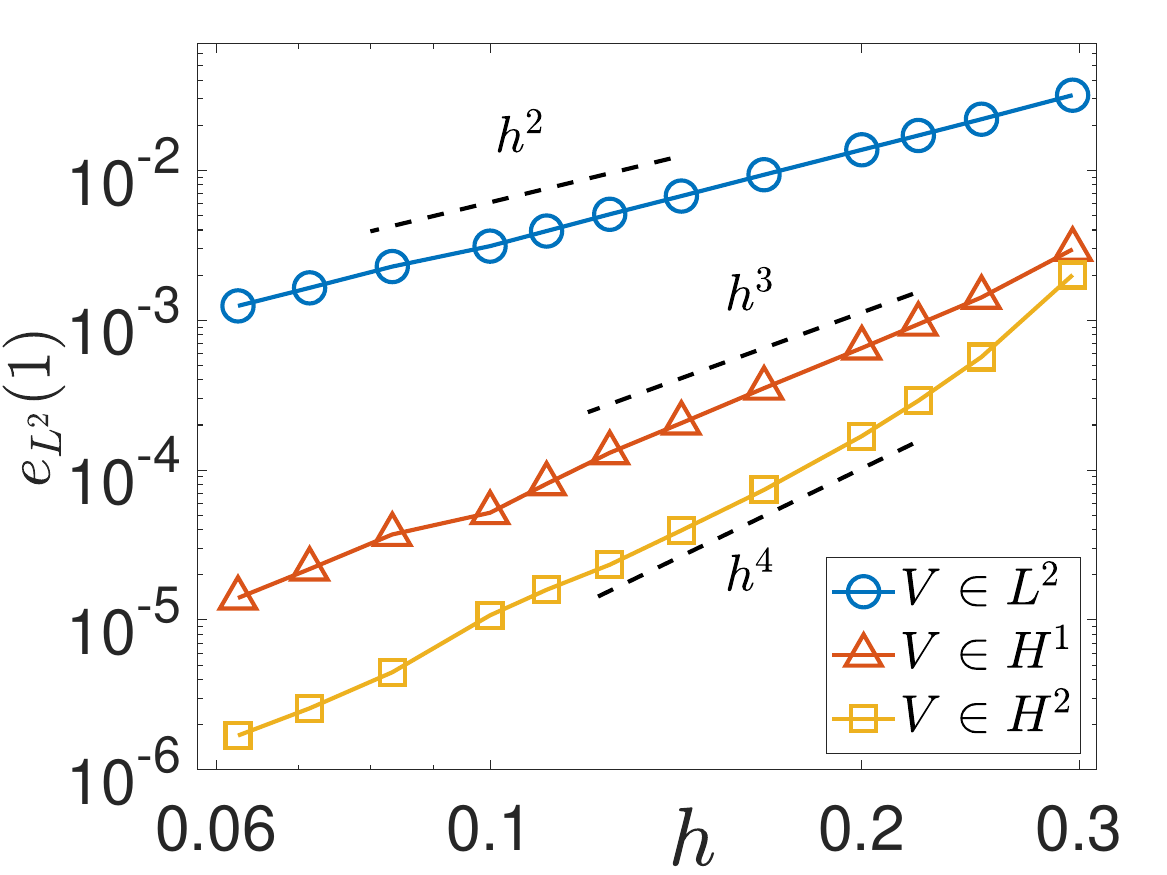}
    \includegraphics[width=0.475\linewidth]{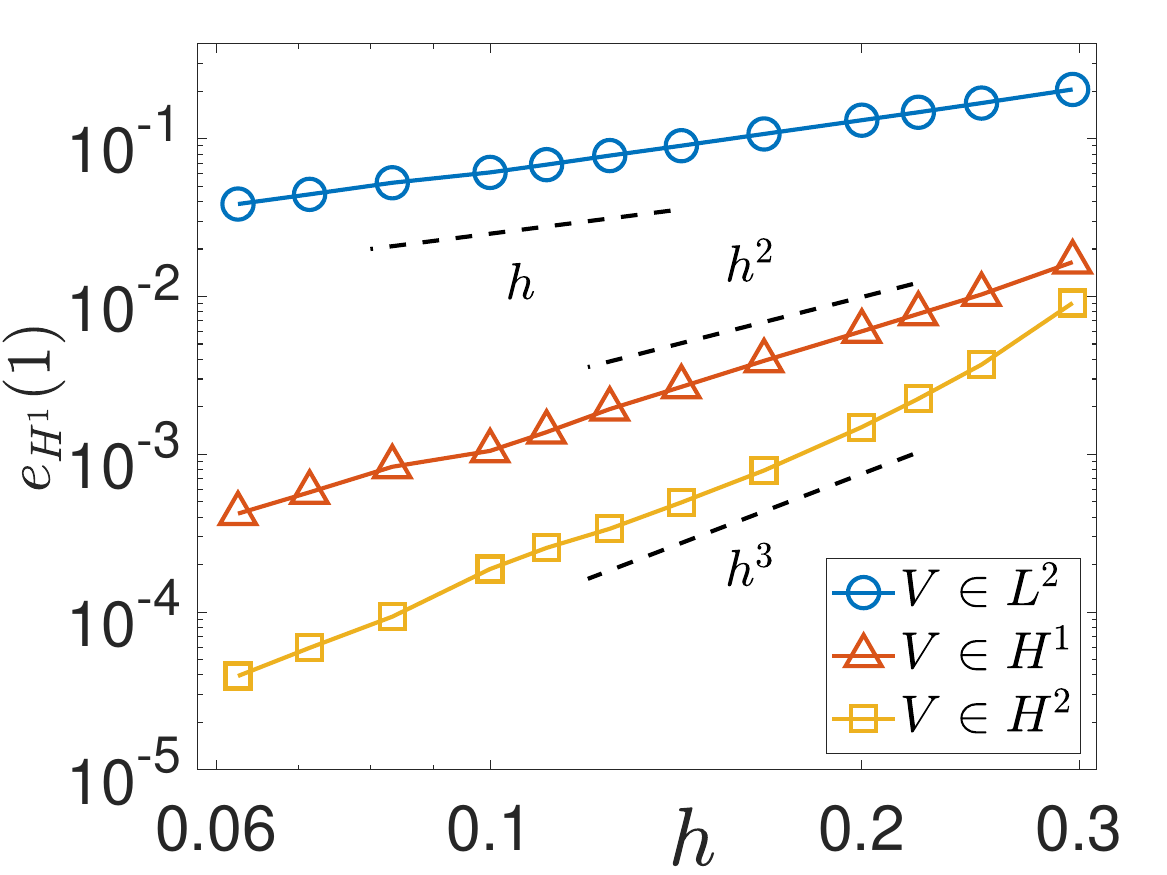}
    \caption{Errors of the PP-FP method for $V=V_{0} \in L^2$, $V=V_{1} \in H^1$ and $V=V_{2} \in H^2$: $L^2$-error (left) and $H^1$-error (right)}
    \label{fig:section12}
\end{figure}

Then we proceed to test the convergence behaviour of the PP-FP method. As the result for $V = V_{1/2}$ is already presented in \cref{fig:section11}, we only consider $V_0$, $V_1$ and $V_2$ in the following. The corresponding numerical results are presented in \cref{fig:section12}. In particular, when $V = V_j \in H^j$, the $L^2$- and $H^1$-norm convergence are $(j+2)$-order and $(j+1)$-order, respectively, for $j=0, 1, 2$.   

According to the results in \cref{fig:section11,fig:section12} for $V$ of different regularity, we see that, in all the cases, the observed convergence orders of the PP-FP method are optimal with respect to the regularity of the underlying exact solutions, even if the potential is low-regularity or singular. In particular, for $H^s$-solutions with $s \geq 2$, the PP-FP method is observed to be $s$-order convergent in the $L^2$-norm and $(s-1)$-order convergent in the $H^1$-norm, independent of the regularity of the potential. Although such a conclusion comes only from the numerical observation, we believe it will hold in general, and a rigorous proof will be presented in our future work. 

Next, we compare the PP-FP method with the FP and FS discretization of the time-splitting methods, using the same setting of $\Omega$, $\beta$, and $\psi_0$ as above. The comparison focuses on accuracy, computational efficiency, and long-time behavior. We use the low-regularity potential $V_{1/2} \in L^\infty$ and the singular potential $V_{0} \in L^2$ (given in \cref{eq::vs1,eq:V0}, respectively) as examples to illustrate the results. 

\begin{figure}
    \centering
    \includegraphics[width=0.475\linewidth]{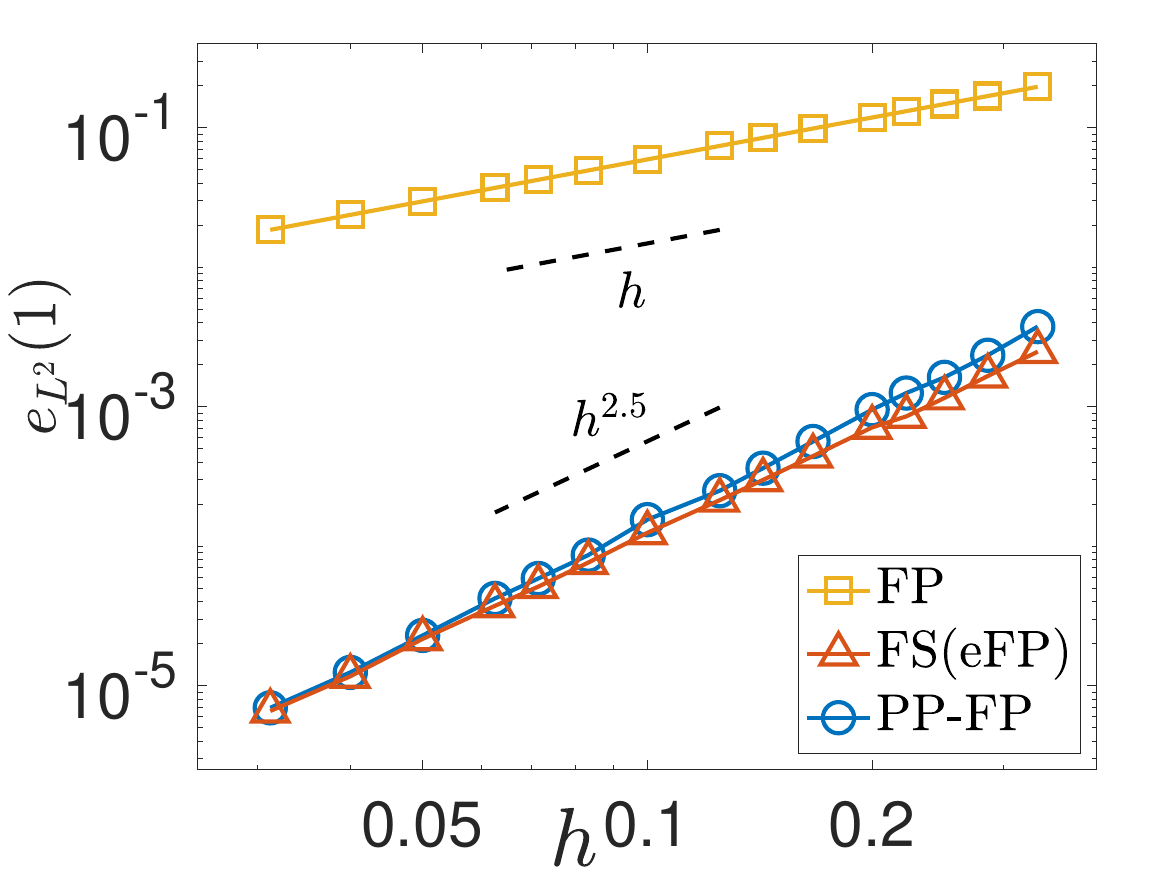}
    \includegraphics[width=0.475\linewidth]{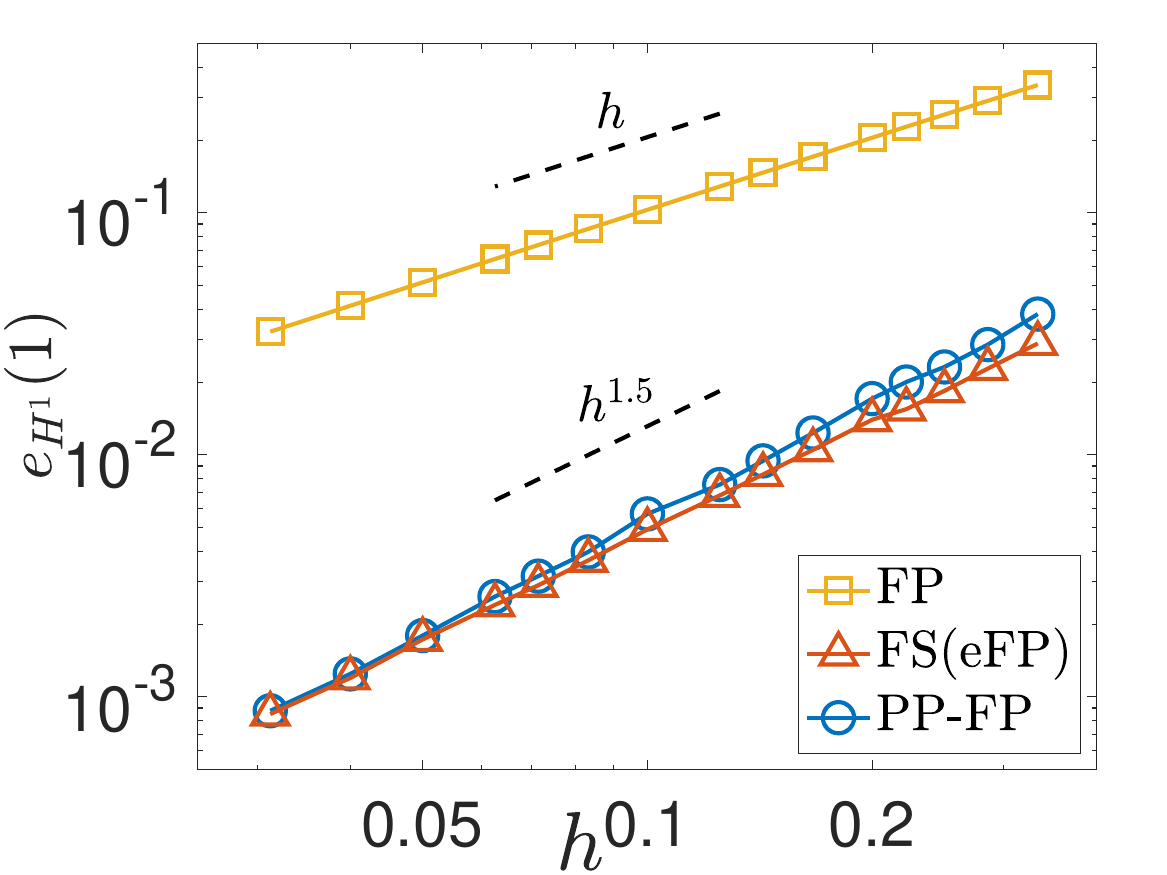}
    \caption{Comparison of the spatial errors of PP-FP, FP, and FS (eFP) methods for the $V = V_{1/2} \in L^\infty$: $L^2$-error (left) and $H^1$-error (right)}
    \label{fig:section42}
\end{figure}

\begin{figure}
    \centering
    \includegraphics[width=0.475\linewidth]{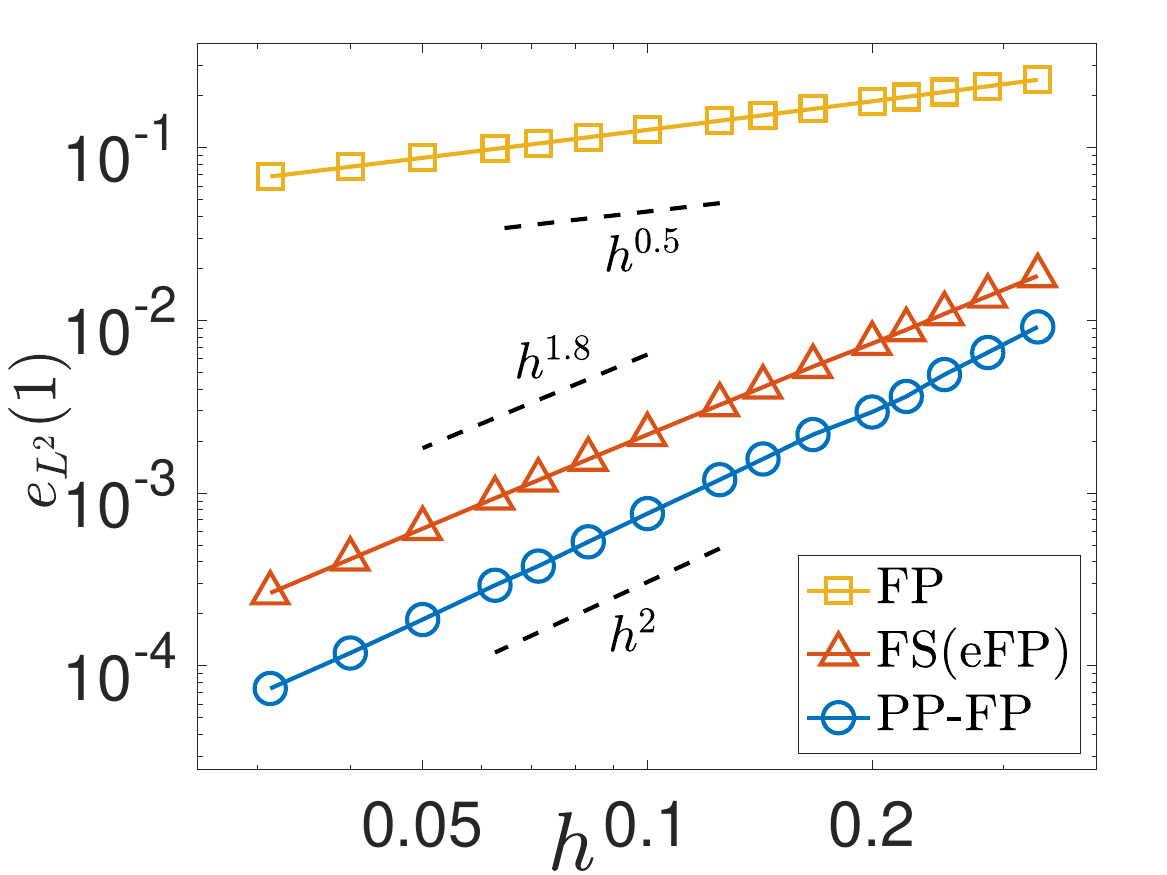}
    \includegraphics[width=0.475\linewidth]{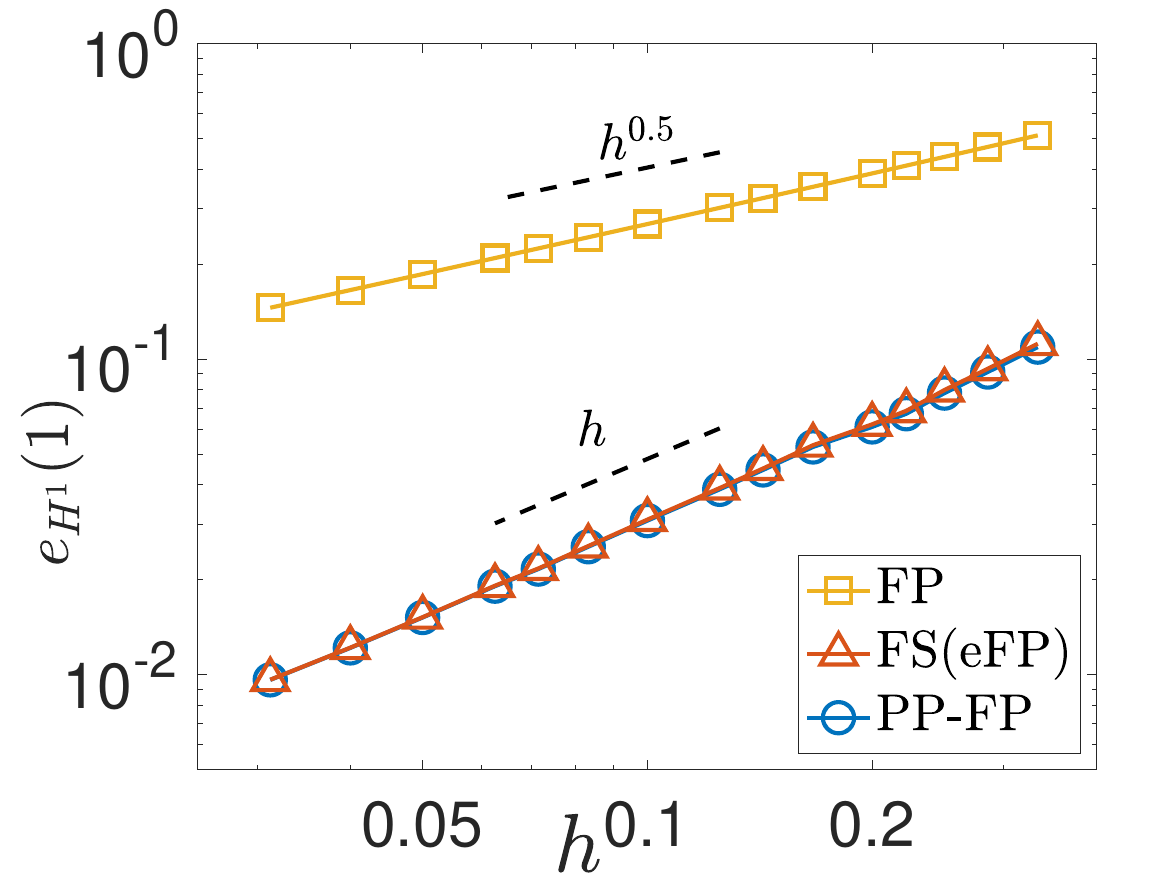}
    \caption{Comparison of the spatial errors of PP-FP, FP, and FS (eFP) methods for the $V = V_0 \in L^2$: $L^2$-error (left) and $H^1$-error (right)}
    \label{fig:section40}
\end{figure}

In the computation, the ST-FP method \cite{bao2003numerical} is used for the FP method, and the ST-eFP method \cite{bao2024efp} is used for the FS method, which is much more efficient than directly evaluating the Fourier integrals. When applying the ST-FP method to the singular potential $V_0$, we shifted the grid by $h/2$ to avoid evaluating $V_0$ at the singular point $x=0$. The reference solution is computed using the ST-PP-FP method with $\tau=\tau_e$ and $h = h_e$. The numerical solutions are computed using the same $\tau = \tau_e$ and varying $h$. The numerical results are exhibited in \cref{fig:section40,fig:section42}. 
% \cw{show results for all the four potentials?}

% Two one-dimensional potentials with different regularity are considered:
% \begin{equation}
% V_s(x)=|x|^{-0.49},
% \quad x\in\Omega=(-16,16),
% \end{equation}
% which is singular and only belongs to $L^2(\Omega)$, and $V_1$ in \cref{eq::vs1} which belongs to $H^1(\Omega)$. In all spatial convergence tests, the time step is chosen as $10^{-6}$, which is sufficiently small so that the temporal discretization error is negligible. The corresponding numerical results are presented in Figs.~\ref{fig:section40} and \ref{fig:section42}. 

% \begin{figure}
%     \centering
%     \includegraphics[width=1\linewidth]{fig//section40.pdf}
%     \caption{Comparison of the spatial errors of PP-FP, FP, and FS (eFP) methods for the $V = V_0 \in L^2$. \cw{no need LT in the legend}}
%     \label{fig:section40}
% \end{figure}

% \begin{figure}
%     \centering
%     \includegraphics[width=1\linewidth]{fig//section42.pdf}
%     \caption{Comparison of the spatial errors of PP-FP, FP, and FS (eFP) methods for the $V = V_{1/2} \in L^\infty$. \cw{no need LT in the legend}}
%     \label{fig:section42}
% \end{figure}

For the low-regularity $L^\infty$-potential, the PP-FP method shows 2.5-order convergence in the $L^2$-norm, and $1.5$-order convergence in the $H^1$-norm, which are the same as the FS (eFP) method, and is much higher than the FP method. 
% Moreover, such convergence order agrees with the (almost) $H^{2.5}$-regularity of the exact solution. 
The same phenomenon is observed also for more regular potentials which are not shown here for brevity. For the singular $L^2$ potential, the difference among the three spatial discretizations is more obvious. The PP-FP method is second-order in the $L^2$-norm and first-order in the $H^1$-norm. The FS (eFP) method shows a slight order reduction to $1.8$ in the $L^2$-norm, while the FP method exhibits only $0.5$-order convergence in the $L^2$- and $H^1$-norms. These indicate that the PP-FP method can achieve (at least) the same accuracy and convergence orders as the FS method, and is significantly more accurate than the FP method. 

In terms of the computational complexity, it is clear that the PP-FP method has exactly the same per-step computational cost as the standard FP method, both of which are of $O(M^d \log M^d)$ in $d$-dimensions. The computational cost of the FS method is more involved and depends on the particular quadrature methods, which is in general much larger than that of the FFT, especially for low-regularity or singular potentials. Even using the proposed eFP method in the implementation of the FS method, the computational cost is at least $3^d$ times the cost of the PP-FP and FP methods as the FFT is performed on an extended grid of length $3M$ in each dimension. Moreover, when combined with the time-splitting methods, the eFP method requires the precomputation of $P_{2M}\left(e^{-i\tau V}\right)$, 
which is much less flexible compared to $P_N V$, as it needs to be recomputed whenever the time step size $\tau$ changes. 

% For the $H^1$ potential, the gap between the three methods becomes smaller; see Fig.~\ref{fig:section42}. Nevertheless, PP-FP remains clearly more accurate than the standard FP method and again gives results comparable to those of eFP.

% These results show that PP-FP consistently outperforms the standard FP discretization and attains the accuracy of eFP, with its advantage being most evident for the rough $L^2$ potential. Importantly, this accuracy is achieved without introducing the enlarged auxiliary grid required by eFP.

% This difference in spatial treatment also leads to a clear distinction in computational cost. Both FP and PP-FP evolve the numerical solution on the original $M^d$-point grid and have a cost of $O(M^d\log M)$ per time step. However, in the implementation for eFP method, the enlarged-grid operations lead to a cost of
% \begin{equation}
% O\left((3M)^d\log(3M)\right)
% \end{equation}
% per time step, with a prefactor that becomes increasingly significant as the spatial dimension grows. 

\begin{figure}
    \centering
    \includegraphics[width=0.475\linewidth]{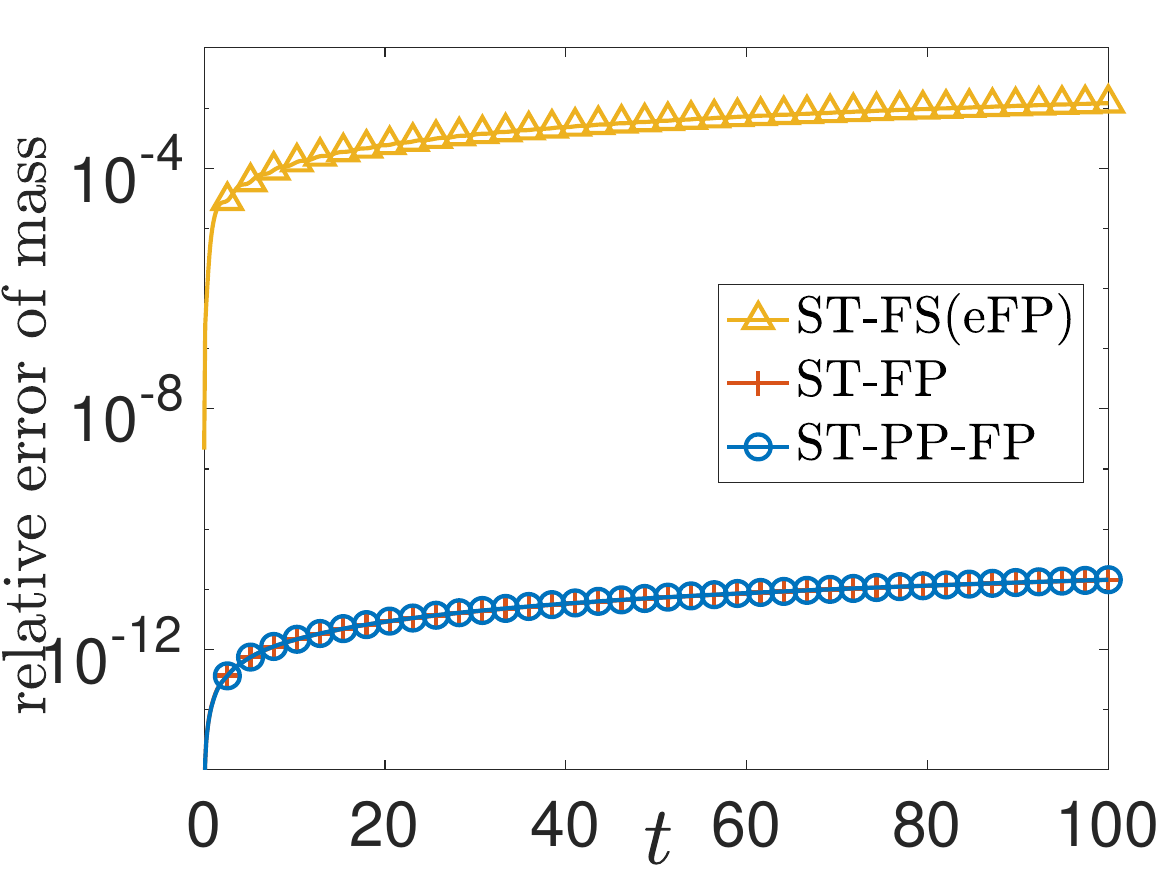}
    \includegraphics[width=0.475\linewidth]{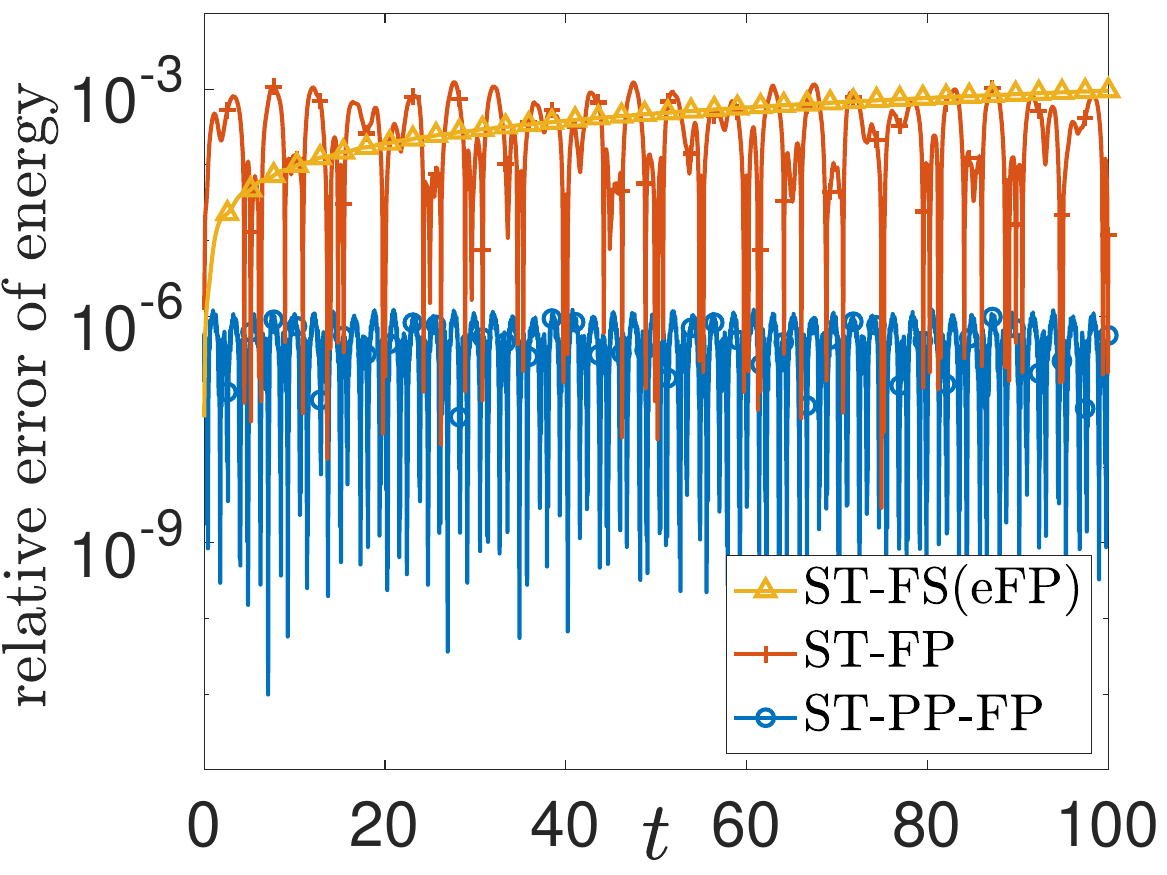}
    \caption{Errors of the mass and energy for the PP-FP, FS, and FP methods with ST for temporal discretization: mass error (left) and energy error (right)}
    \label{fig:mass_energy_1d_2x2}
\end{figure}

Beyond accuracy and computational efficiency, the PP-FP method also exhibits excellent performance in the long-time simulation due to its structure-preserving property. In the following, we compare the error of the mass and energy of the numerical solutions over a long time interval $t \in [0, 100]$. 
% obtained by the ST-PP-FP (for PP-FP), ST-FP (for FP), and ST-eFP (for FS) methods. 
We choose the $L^\infty$-potential $V_{1/2}$ as an example, and use the same $\tau = 10^{-3}$ and $h = 2^{-4}$ for all the three methods. The potential energy is computed by evaluating the corresponding integral exactly for the three methods (which is important due to the low regularity of the potential). This can be done by using Parseval's identity. 
The corresponding results are shown in \cref{fig:mass_energy_1d_2x2}. 
% \cw{The difference between FP and PP-FP should be much more obvisous, and should be smaller.}

% \begin{figure}
%     \centering
%     \includegraphics[width=1\linewidth]{fig//mass_energy_1d_2x2.pdf}
%     \caption{Errors of the mass and enenrgy for the PP-FP, FP, and FS methods with ST for temporal discretization \cw{energy error: three methods in one figure, one $\tau$ is enough, no need to put $\tau$ in legend}}
%     \label{fig:mass_energy_1d_2x2}
% \end{figure}

We observe that both PP-FP and FP preserve the mass up to machine precision but FS (eFP) does not. For the energy, PP-FP produces a small and oscillating error due to its improved spatial accuracy and the preservation of the time symmetry of the temporal integrator. In comparison, the energy error of the FP method, though it also oscillates without growing in time due to the time symmetry, the value is about three magnitudes larger than that of the PP-FP. The FS (eFP) produces an energy error smaller than the FP initially but it visibly grows in time because the time symmetry is lost. 

% The results clearly indicate a structural difference between the three methods at the level of mass conservation. Both FP and PP-FP preserve the discrete mass exactly up to machine precision throughout the computation. In contrast, the eFP method does not preserve the mass exactly.

% For the discrete energy, none of the three methods is exactly energy-conserving, as expected for a first-order splitting scheme. Nevertheless, their long-time energy behaviors are visibly different. For the eFP method, the energy error displays a persistent increasing trend, indicating an accumulation of energy drift over long times. By contrast, the energy errors of FP and PP-FP remain bounded and oscillatory throughout the computation. Moreover, the magnitude of the energy error produced by PP-FP is clearly smaller than that of the standard FP method in the reported tests.

% These observations show that PP-FP retains the exact mass conservation property of the standard FP method while improving the long-time behavior of the discrete energy. Compared with eFP, it avoids the visible energy drift; compared with FP, it yields a smaller and more stable energy error. Therefore, from the viewpoint of conservation and long-time stability, PP-FP preserves the main structural advantage of FP and at the same time provides a more favorable treatment of low-regularity potentials.

In summary, compared with FS (eFP), PP-FP attains comparable short-time accuracy at a substantially lower computational cost and exhibits better long-time behavior in the reported tests. Compared with FP, PP-FP is more accurate in both short- and long-time simulations at the same per-step cost. These observations are summarized in Table~\ref{tab:fp_efp_fpp_summary}.

\begin{table}
\centering
\caption{Comparisons of PP-FP, FS, and FP methods combined with (symmetric) time-splitting methods}
\label{tab:fp_efp_fpp_summary}
\begin{threeparttable}
\small
\setlength{\tabcolsep}{4.5pt}
\renewcommand{\arraystretch}{1.15}
\begin{tabular}{@{}
>{\raggedright\arraybackslash}p{0.22\textwidth}
>{\raggedright\arraybackslash}p{0.15\textwidth}
>{\centering\arraybackslash}p{0.15\textwidth}
>{\centering\arraybackslash}p{0.18\textwidth}
>{\centering\arraybackslash}p{0.15\textwidth}@{}}
\toprule
\multicolumn{2}{@{}l}{}
& \textbf{PP--FP} & \textbf{FS(eFP)} & \textbf{FP} \\
\midrule
\multirow{2}{0.22\textwidth}{Convergence: $V \in L^2$}
& $L^2$-order & $2$ & $1.8$ & 0.5 \\
\cmidrule(lr){2-5}
& $H^1$-order & $1$ & $1$ & 0.5 \\
\midrule
\multirow{2}{0.22\textwidth}{Convergence: $V \in L^\infty\cap H^{\frac12-}$}
& $L^2$-order & $2.5$ & $2.5$ & $1$ \\
\cmidrule(lr){2-5}
& $H^1$-order & $1.5$ & $1.5$ & $1$ \\
\midrule
\multicolumn{2}{@{}l}{Cost per time step}
& $O(M^d\log M)$ & $O((3M)^d\log M)$ & $O(M^d\log M)$ \\
\addlinespace[2pt]
\multicolumn{2}{@{}l}{Mass conservation}
& Yes & No & Yes \\
\addlinespace[2pt]
\multicolumn{2}{@{}l}{Energy conservation}
& No & No & No \\
\addlinespace[2pt]
\multicolumn{2}{@{}l}{Time symmetry}
& Yes & No & Yes \\
\bottomrule
\end{tabular}

% \begin{tablenotes}[flushleft]
% \footnotesize
% \item A dash indicates that no clear convergence order is observed.
% \end{tablenotes}

\end{threeparttable}
\end{table}

\subsection{Temporal convergence}\label{sec4:splitting}

We study the temporal convergence of the TS-PP-FP method, and systematically investigate how the performance of the time-splitting methods depends on the regularity of the potential. We aim at answering the following two questions: (i) for a potential of given regularity, what is the convergence order of the time-splitting methods and (ii) for a time-splitting method to achieve a given convergence order, what is the minimal regularity required for the potential. We fix $\Omega=(-16,16)$, $T=1$, and consider the linear case with $\beta = 0$ in \cref{eq::nlse}. The same conclusion can be drawn for the nonlinear case. Again, the initial datum is chosen as a Gaussian $\psi_0(x)=e^{-x^2/2}$ for all the tests.  
% The reference solution is computed using the time-splitting PP-FP method with $\tau=\tau_e=10^{-6}$ and $h=h_e=2^{-9}$. 

\begin{figure}
    \centering
    \includegraphics[width=0.475\linewidth]{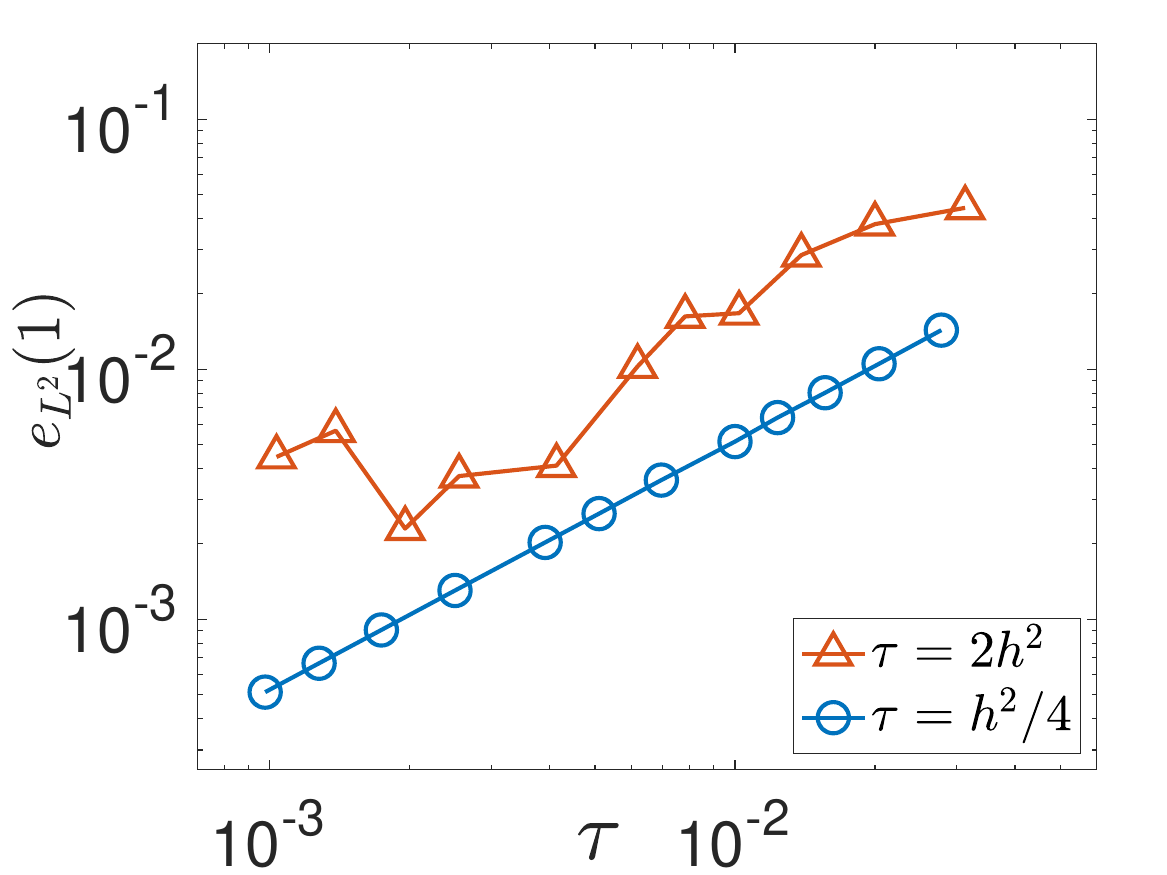}
    \includegraphics[width=0.475\linewidth]{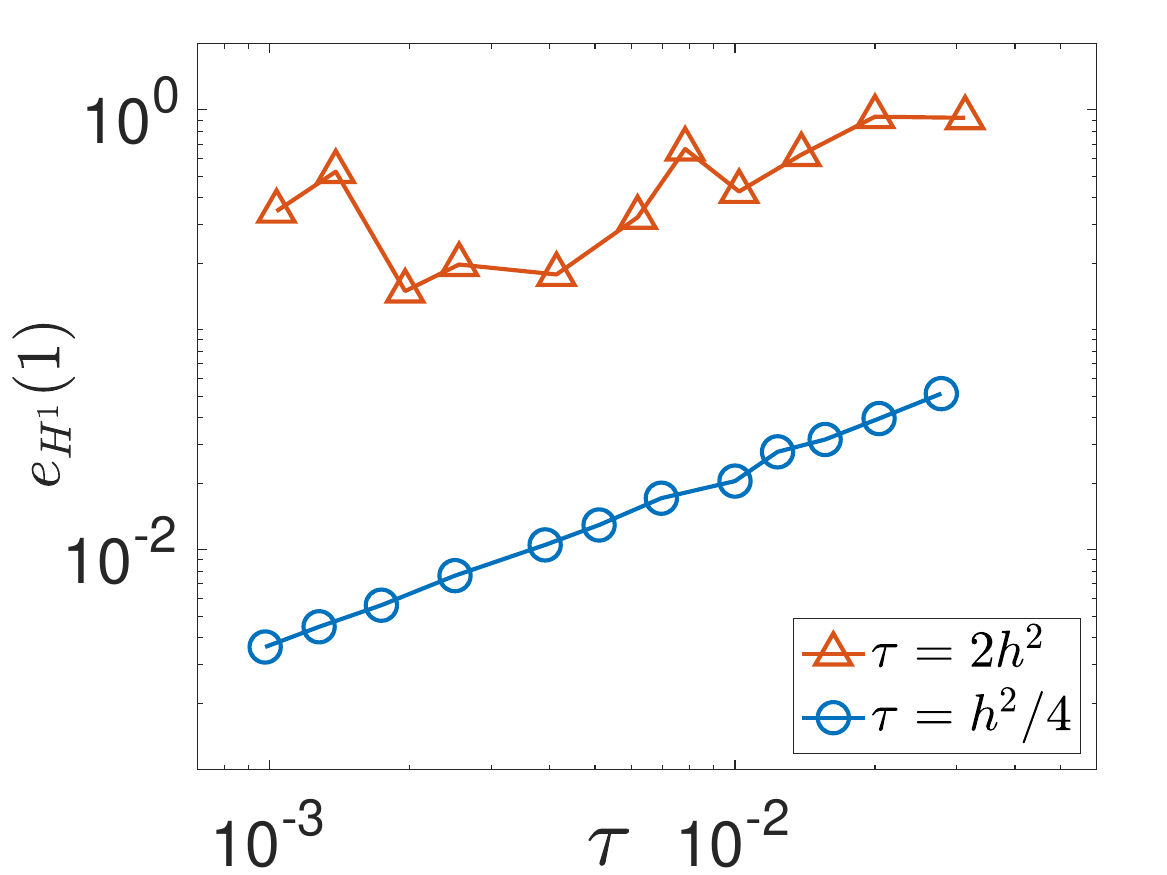}
    \caption{Comparison of the errors of LT-PP-FP method with $\tau = h^2/4$ and $\tau = 2h^2$: $L^2$-error (left) and $H^1$-error (right)}
    \label{fig:section60}
\end{figure}

We first present a numerical experiment to show that a CFL-type time step size restriction $\tau \leq h^2/\pi$ (which is required in the error estimates for the time-splitting FS and eFP methods in \cite{bao2024efp,bao2024tsnls}) is also needed when implementing the TS-PP-FP methods. We use the $L^\infty$-potential $V_{1/2}$ in \cref{eq::vs1} to illustrate this. We plot the errors for the LT-PP-FP method with $\tau = h^2/4$ and $\tau = 2h^2$ in \cref{fig:section60}. We see that when the CFL-type condition is violated, the behaviour of the LT-PP-FP method is very irregular and the values of the errors are much larger even though a finer grid is used. In the following, we will always choose $\tau$ and $h$ such that $\tau = h^2/4$.
% the regularity requirement on the potential can be significantly reduced, or, equivalently, the convergence order can be significantly improved (if not yet optimal). Interestingly, such improvement is not possible when using the FP method for the spatial discretization, indicating that different choices of the spatial discretization can also have influence on the bahaviour of time-spllting methods. 
% \cw{Show the numerical results to compare the temporal errors for TS-PP-FP with $\tau \leq h^2$, TS-PP-FP (with $\tau \gg h^2$), TS-FP with $\tau \leq h^2$, TS-FP (with $\tau \gg h^2$). Two figures: $L^2$ and $H^1$}

\begin{figure}
    \centering
    \includegraphics[width=0.475\linewidth]{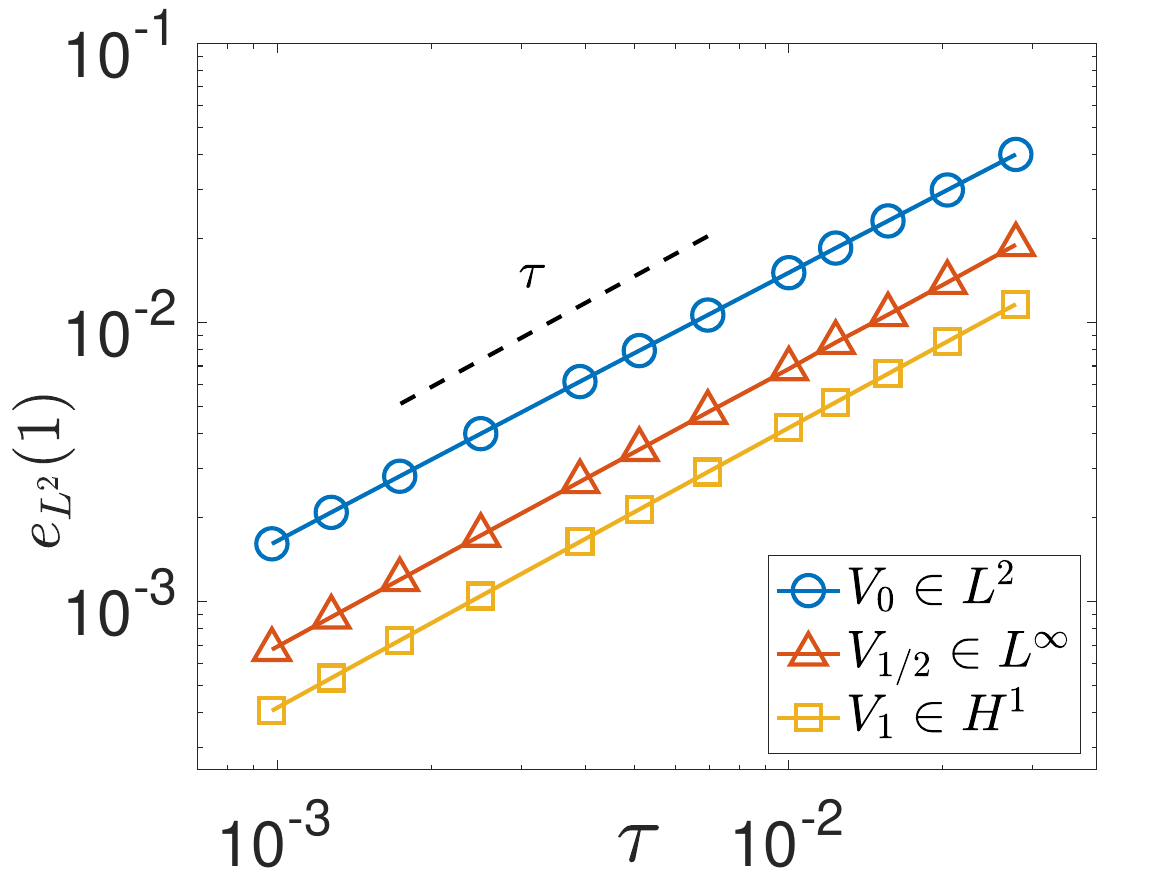}
    \includegraphics[width=0.475\linewidth]{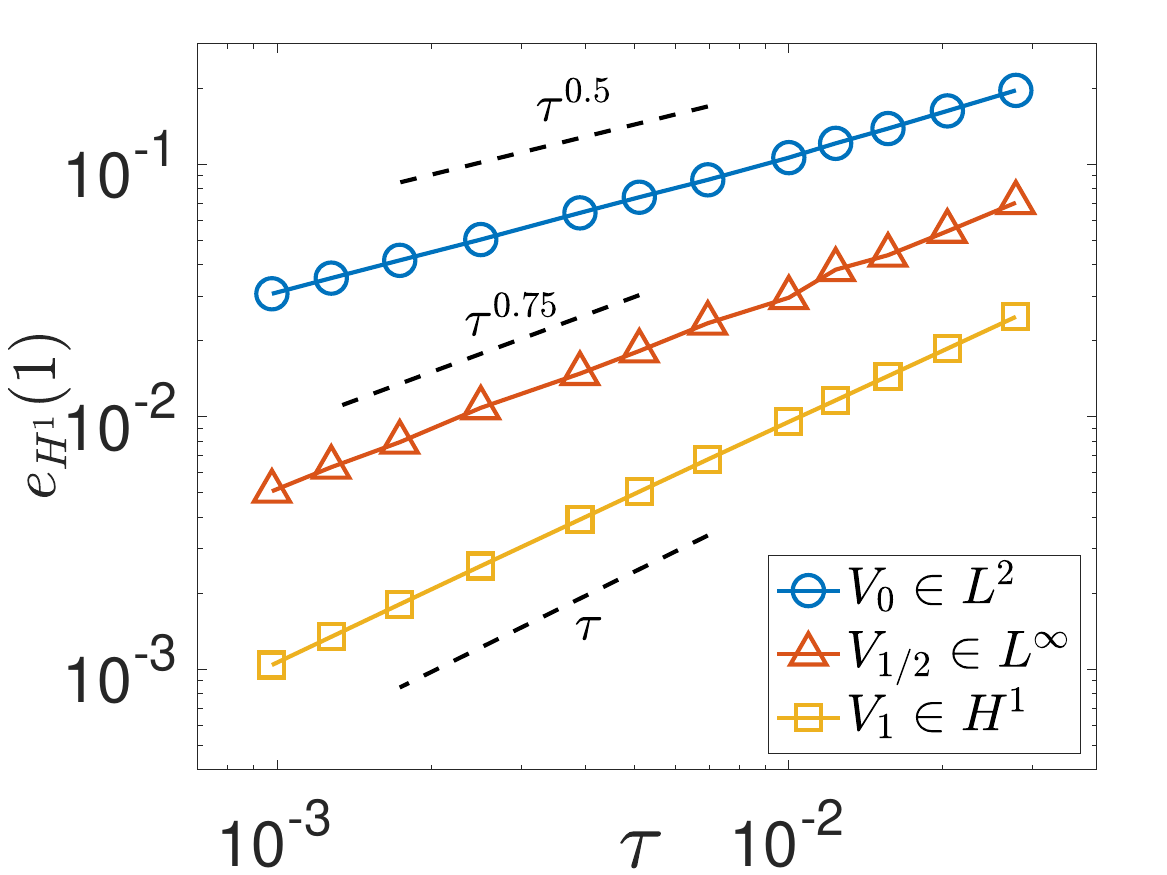}
    \caption{Temporal errors of the LT-PP-FP method for potentials with different regularity: $L^2$-error (left) and $H^1$-error (right)}
    \label{fig:section61}
\end{figure}

\begin{figure}
    \centering
    \includegraphics[width=0.475\linewidth]{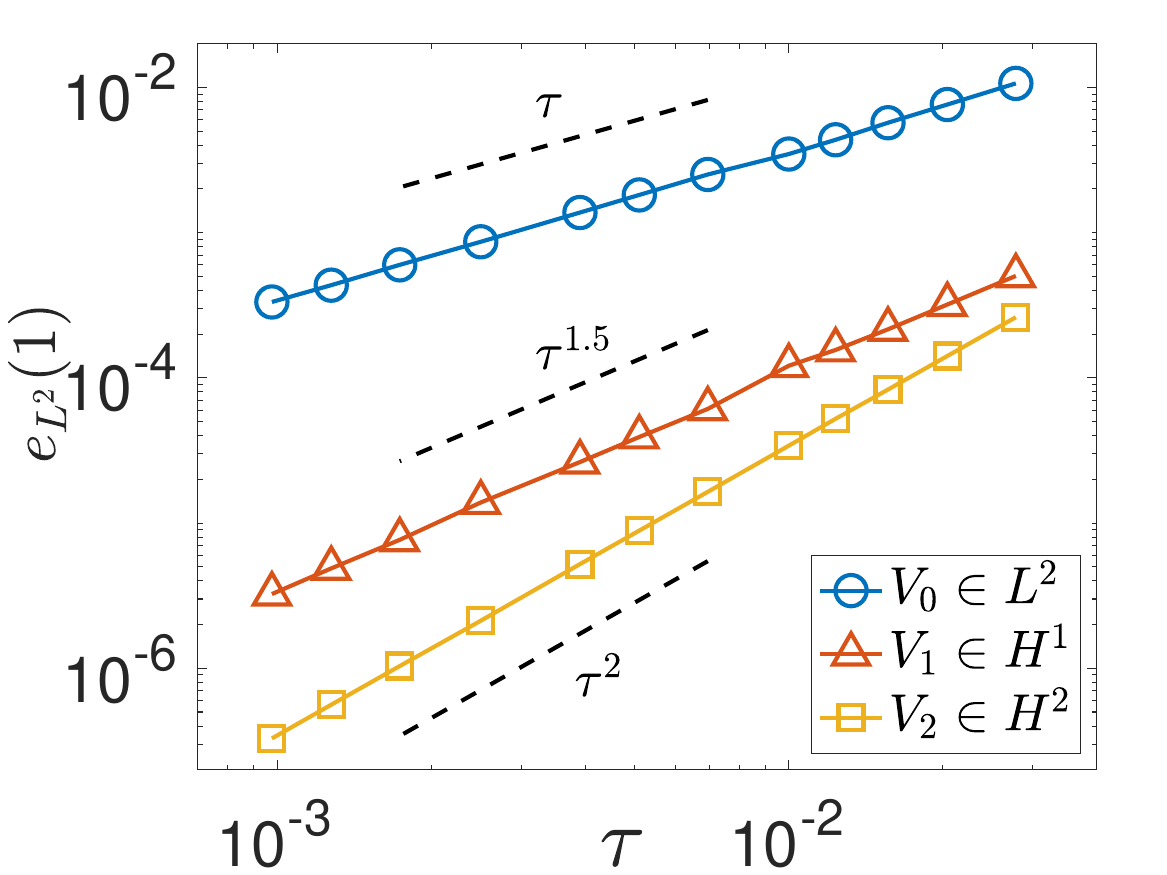}
    \includegraphics[width=0.475\linewidth]{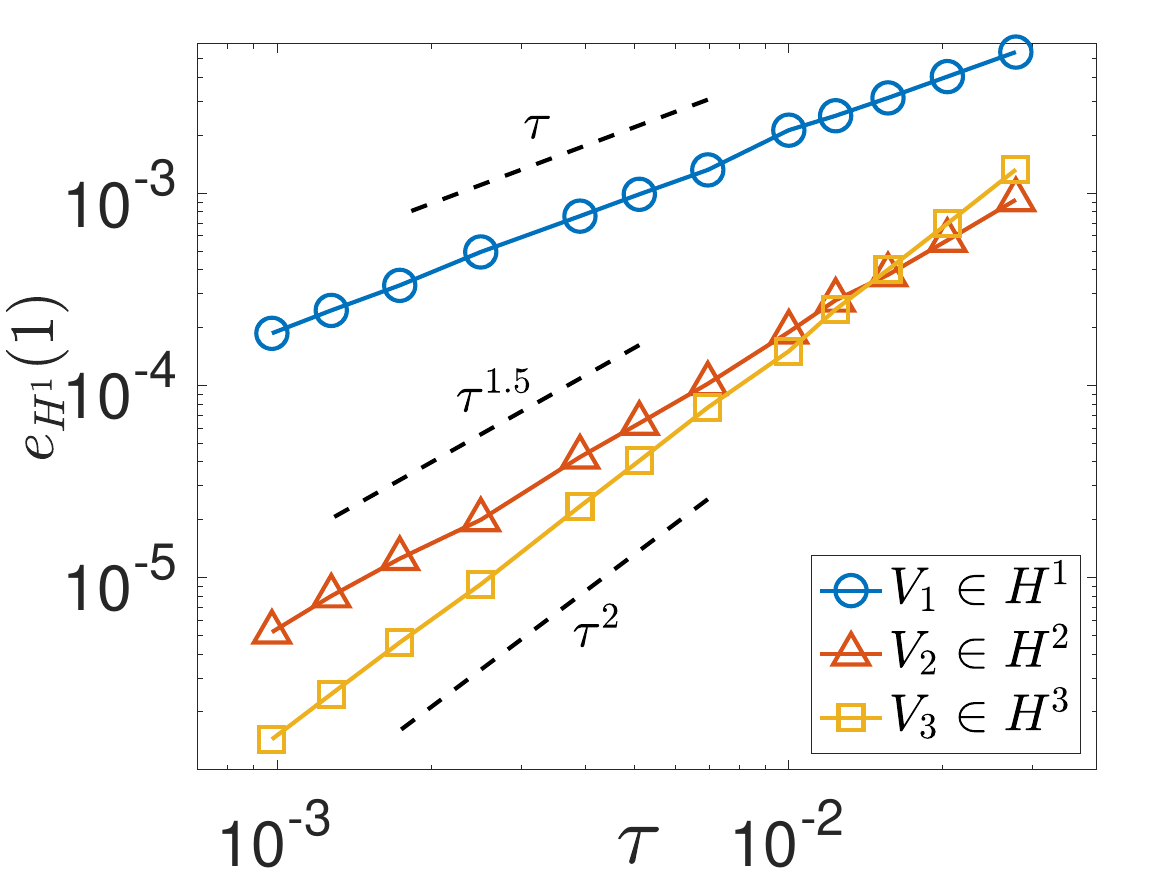}
    \caption{Temporal errors of the ST-PP-FP method for potentials with different regularity: $L^2$-error (left) and $H^1$-error (right)}
    \label{fig:section62}
\end{figure}

\begin{figure}
    \centering
    \includegraphics[width=0.475\linewidth]{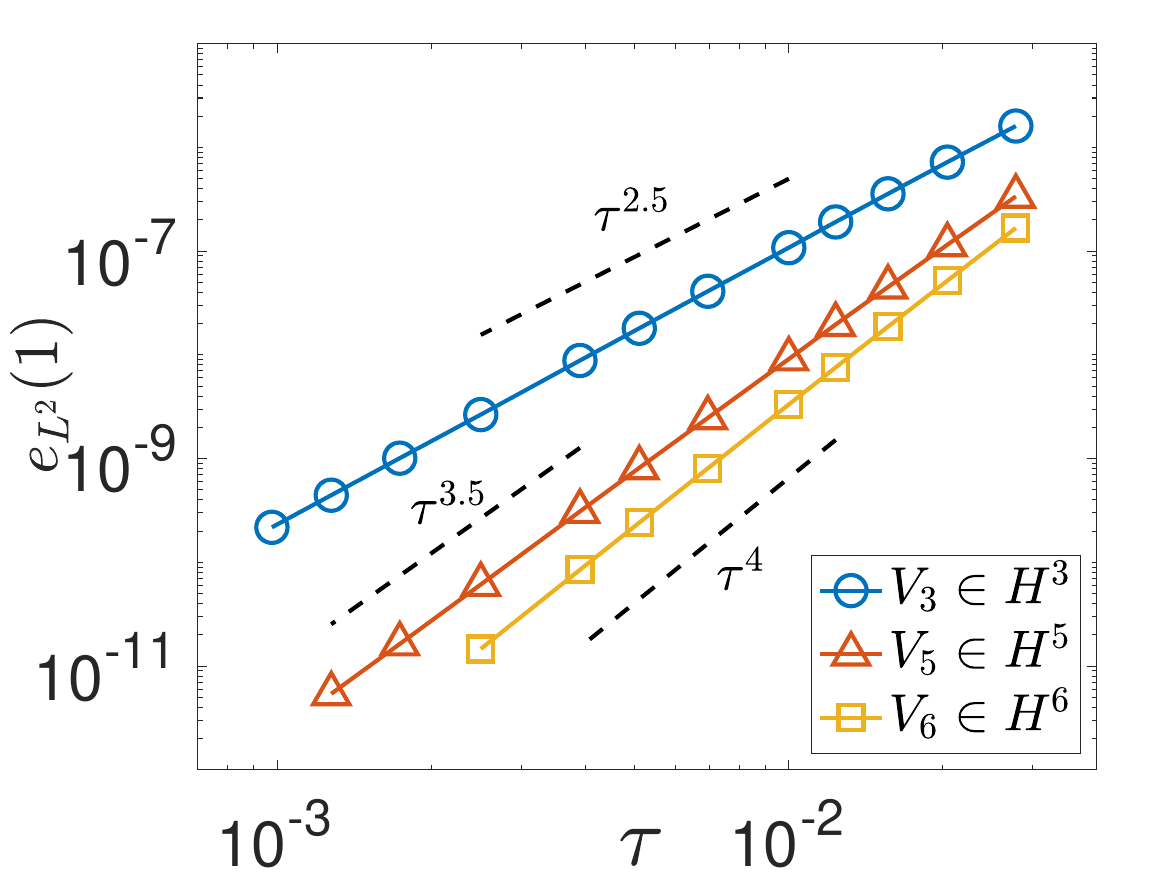}
    \includegraphics[width=0.475\linewidth]{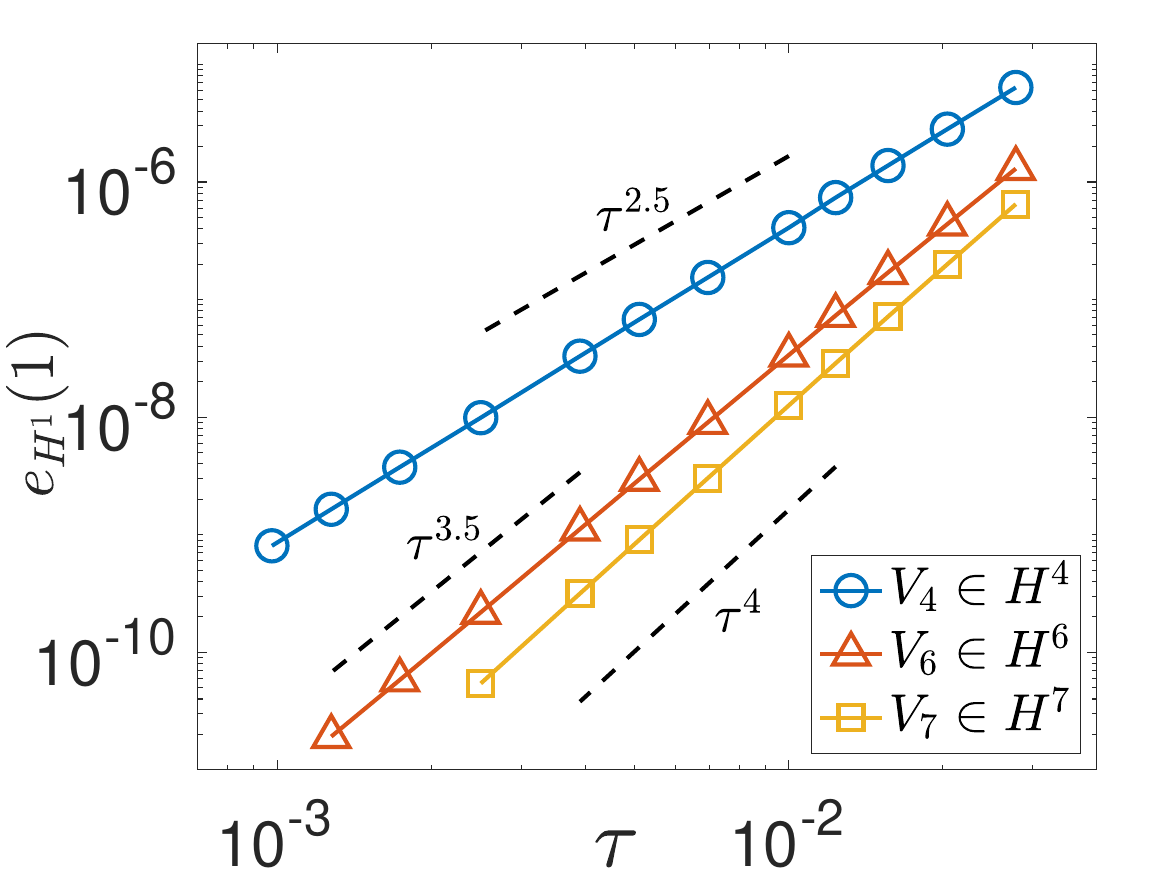}
    \caption{Temporal errors of the S4-PP-FP method for potentials with different regularity: $L^2$-error (left) and $H^1$-error (right)}
    \label{fig:section63}
\end{figure}

\begin{figure}
    \centering
    \includegraphics[width=0.475\linewidth]{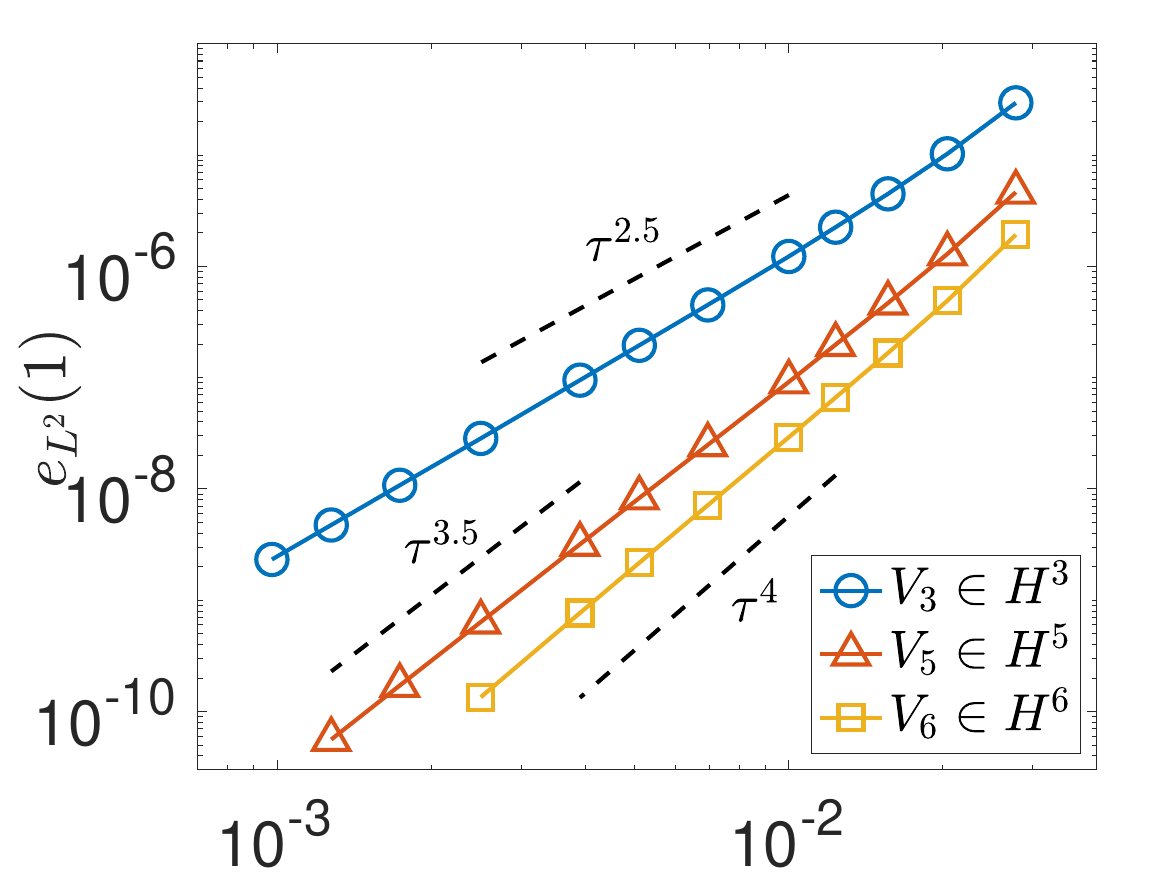}
    \includegraphics[width=0.475\linewidth]{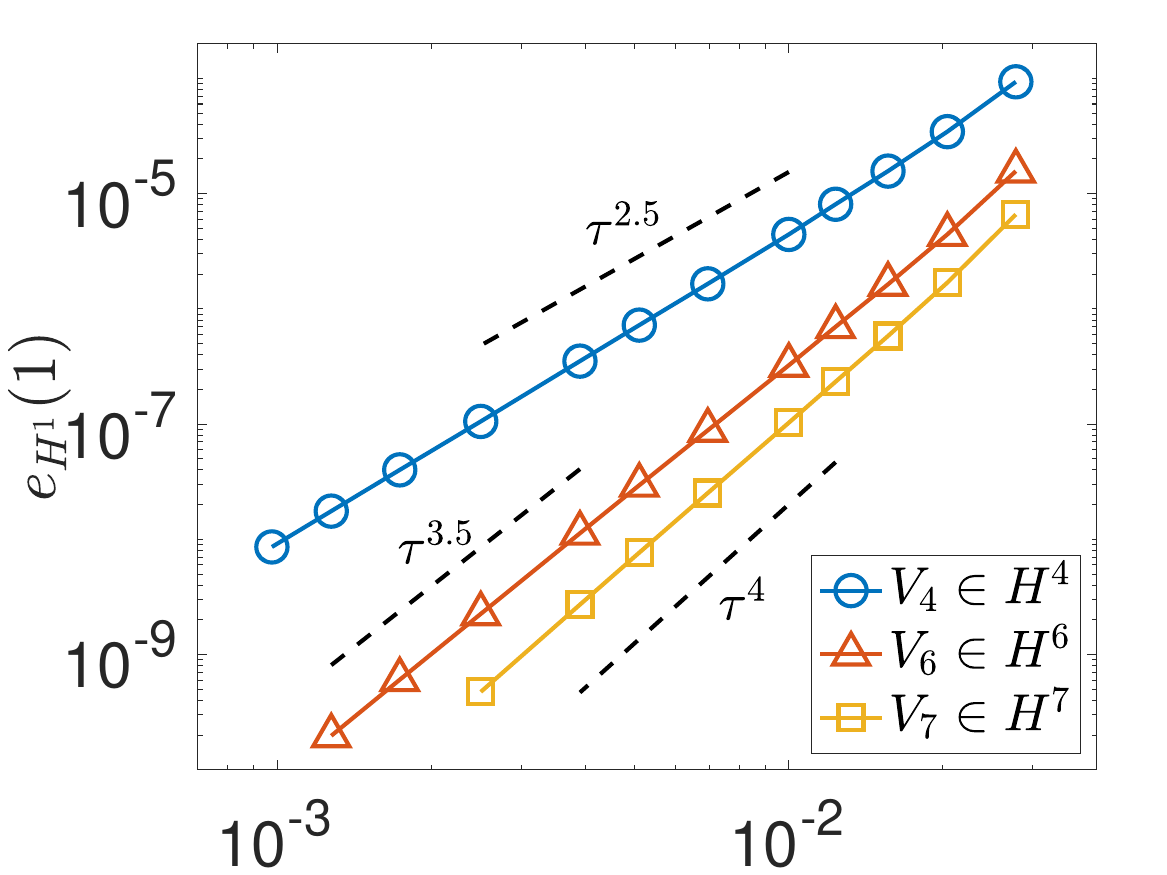}
    \caption{Temporal errors of the S4c-PP-FP method for potentials with different regularity: $L^2$-error (left) and $H^1$-error (right)}
    \label{fig:section64}
\end{figure}

We will consider the potentials $V = V_j$ with $j=0, 1/2, 1, \cdots, 7$, where $V_j$ for $j=0, 1/2, 1, 2$ are given in \cref{eq::vs1,eq:V0}, and, for $j = 3, \cdots, 7$, are given in terms of their Fourier coefficients as
\begin{equation}
    \widehat{(V_j)}_l = A (-1)^l \frac{\langle l/l_0 \rangle^{-(\frac{1}{2}+j)}}{\ln{(2+|l/l_0|)}}, \qquad l \in \mathbb{Z}, 
\end{equation}
where $\langle l \rangle := (1+|l|^2)^{1/2}$, $A>0$ is a normalization constant, and $l_0 = 16$. 
% \begin{equation}
% \label{eq::vs_high}
% \begin{aligned}
% % V_2(x) = |x|^{1.51}\left(1 - \frac{x^2}{16^2}\right)^2, 
% % \quad 
% V_3(x) &= |x|^{2.51}\left(1 - \frac{x^2}{16^2}\right)^3,
% \quad V_4(x)=|x|^{3.51}\left(1 - \frac{x^2}{16^2}\right)^4,\\
% V_5(x) &= \left|\frac{x}{16}\right|^{1.5}|x|^{3.01}
% \left(1 - \frac{x^2}{16^2}\right)^5,\\
% V_6(x) &= \left|\frac{x}{16}\right|^{2}|x|^{3.51}
% \left(1 - \frac{x^2}{16^2}\right)^6,
% \quad V_7(x)=\left|\frac{x}{16}\right|^{3}|x|^{3.51}
% \left(1 - \frac{x^2}{16^2}\right)^7.
% \end{aligned}
% \end{equation}
It can be checked that all the $V_j$ are real-valued and
\[
    V_j\in H^j_{\mathrm{per}}(\Omega),\quad j=0, 1, 2,\ldots,7, \quad V_{1/2}\in L^\infty(\Omega) \cap H^{\frac{1}{2}-}_{\mathrm{per}}(\Omega). 
\]
% \cw{figure for LT $V \in H^{0, 1/2, 1}$ in $L^2$ and $H^1$}
The errors in $L^2$- and $H^1$-norms of the four time-splitting PP-FP methods with potential $V$ of different regularity are presented in \cref{fig:section61,fig:section62,fig:section63,fig:section64}. 
% The numerical observations are summarized in \cref{tab:temporal_orders_fppp}. 

% We then turn to the fourth-order partitioned Runge-Kutta splitting scheme 4RK-PP-FP. The numerical results in Fig.~\ref{fig:section62} show a substantially stronger regularity requirement, which is natural for a higher-order method. More precisely, the $L^2$ error exhibits fourth-order convergence in time when $V\in H^5$, while the $H^1$ error reaches fourth-order convergence when $V\in H^6$. Hence, for this fourth-order splitting strategy, the potential projection treatment still allows the expected high-order temporal behavior to be observed, but only when the regularity of the potential is sufficiently high. At the same time, these threshold values coincide with those previously observed for the extended Fourier pseudo-spectral strategy. This shows that, also for higher-order splittings, PP-FP is able to reduce the regularity requirement to essentially the same level as eFP, while avoiding the oversampling cost of the latter.

% We also present a numerical illustration for the fourth-order compact splitting method CS-PP-FP in Fig.~\ref{fig:section63}.

From the numerical results, with sufficiently smooth initial data, under the CFL-type condition $\tau = h^2/4$,
we can conclude the following:

\medskip 
(i)  When $V \in H^s_\text{per}$ with $s \geq 0$, then the exact solution admits sharp $H^{s+2}_\text{per}$ regularity. The convergence order of the $p$-th order time-splitting PP-FP (TS-PP-FP) method is
\begin{equation}\label{eq:TSPPFP-1}
    \min\{p, 1+s/2\} \text{ in the $L^2$-norm}, \quad \min\{p, 1/2+s/2\} \text{ in the $H^1$-norm}.
\end{equation}
Thus \cref{eq:TSPPFP-1} indicates that TS-PP-FP can always achieve optimal temporal convergence with respect to the regularity of the exact solution. 
On the other hand, under the same setting, the minimal regularity of the potential that is needed for the $p$-th order TS-PP-FP method to achieve optimal order is 
\begin{equation}\label{eq:TSPPFP-223}
V\in H_{\mathrm{per}}^{2p-2}\quad\text{for the }L^2\text{-error},\qquad V\in H_{\mathrm{per}}^{2p-1}\quad\text{for the }H^1\text{-error}.
\end{equation}

(ii)  The above conclusions should also apply to the TS-FS (eFP) methods, possibly excluding the singular cases $V \not \in L^\infty$. 
    
(iii) However, for the TS-FP methods, when $V \in H^s_\text{per}$ with $s \geq 2$, the convergence order of the $p$-th order TS-FP method is
\begin{equation}\label{eq:TSPPFP-2}
    \min\{p, s/2\} \text{ in the $L^2$-norm}, \quad \min\{p, -1/2+s/2\} \text{ in the $H^1$-norm}.
\end{equation}
To achive the same order or accuracy, it needs two additional orders of regularity on $V$, i.e., one needs $V \in H^{s+2}_\text{per}$ for \cref{eq:TSPPFP-1} to hold, or, to obtain optimal convergence orders for the $p$-th order splitting, one needs   
\begin{equation}\label{eq:TSPPFP-224}
V\in H_{\mathrm{per}}^{2p}\quad\text{for the }L^2\text{-error},\qquad V\in H_{\mathrm{per}}^{2p+1}\quad\text{for the }H^1\text{-error}.
\end{equation}

\subsection{Application to the dynamics under Coulomb potentials}

We finally present a three-dimensional example involving a singular multi-center Coulomb potential to demonstrate that the time-splitting PP-FP method can stably simulate the dynamics modelled by the NLSE with a singular Coulomb potential in 3D. 
% Such a setting is representative of molecular-type models, where the electron density may repeatedly concentrate near different nuclei and transfer between them during the evolution.

We consider the NLSE \cref{eq::nlse} with $\beta = -1$ on a domain $\Omega=(-8,8)^3$. The potential is given by the following three-centered Coulomb potential
\begin{equation*}
% \label{eq::multi_coulomb_potential}
    V(\mathbf{x})
    =
    -5\sum_{m=1}^{3}
    \frac{1}{|\mathbf{x}-\mathbf{R}_m|},
    \quad
    \mathbf{R}_1=(-1,-1,0)^T,\
    \mathbf{R}_2=(1,1,0)^T,\
    \mathbf{R}_3=(0,0,0)^T.
\end{equation*}
% Since \eqref{eq::multi_coulomb_potential} is singular at the nuclear centers, this example goes beyond the bounded-potential assumptions used in the theoretical analysis. It is therefore intended as a numerical exploration for the potential projection strategy in a genuinely singular three-dimensional regime.
The initial wave function is chosen as a localized Gaussian packet placed to the left of the three centers,
\begin{equation*}
% \label{eq::multi_coulomb_initial}
    \psi_0(\mathbf{x})=A\exp\left(-\frac{|\mathbf{x}-\mathbf{x}_0|^2}{2s^2}
    +ik_0(x_1-x_{0,1})\right),\quad
    \mathbf{x}_0=(-2,0,0)^T,\ s=0.4,\ k_0=1,
\end{equation*}
where $x_{0,1}$ is the first component of $\mathbf{x}_0$ and $A$ is the amplitude used in the computation. The initial momentum drives the packet toward the Coulomb centers. During the evolution, the attractive singular centers induce deformation, splitting, trapping, and transfer of the density profile between different nuclei regions.

\begin{figure}
    \centering
    \includegraphics[width=1\linewidth]{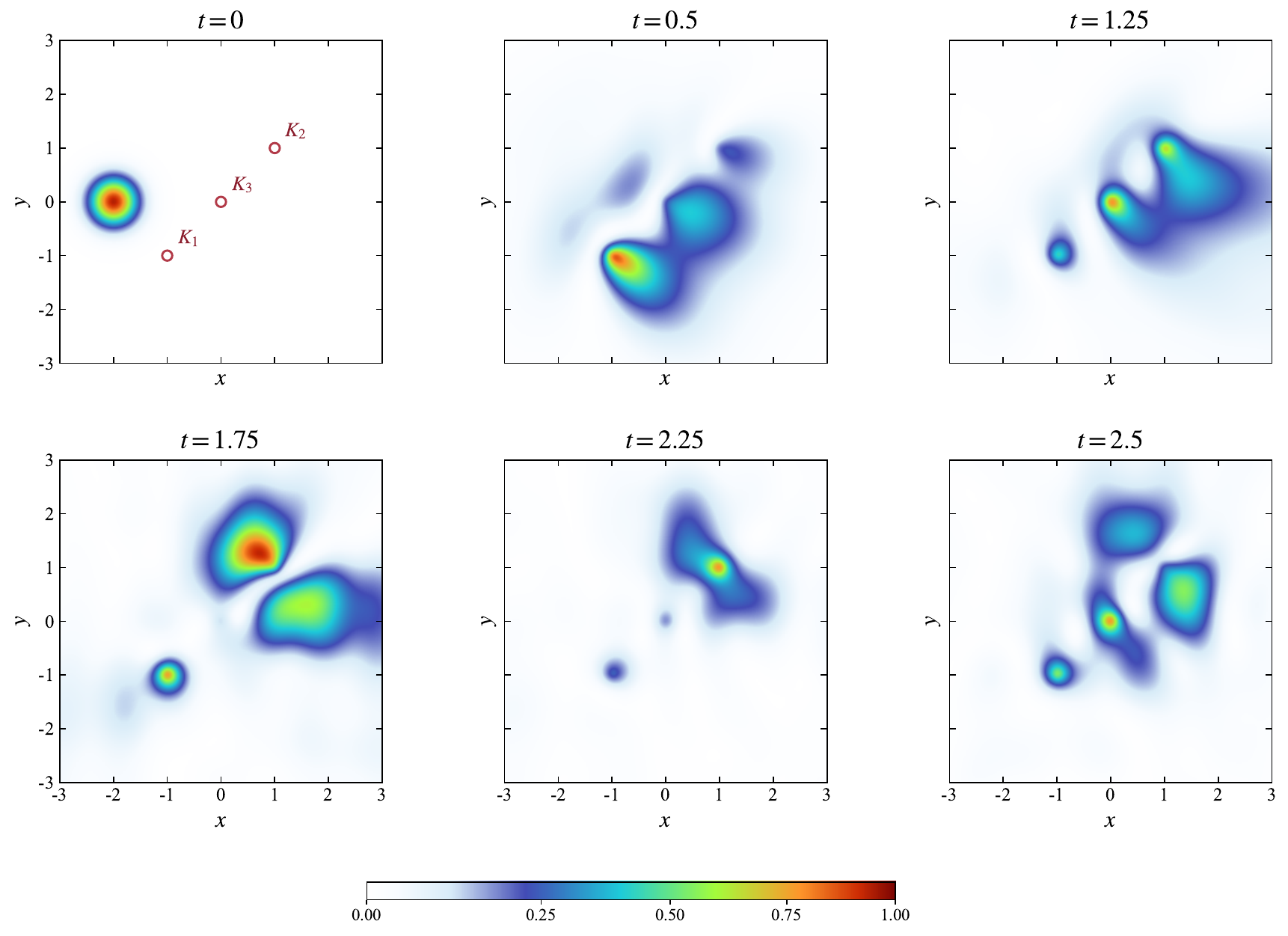}
    \caption{Snapshots of the density $|\psi(\mathbf{x},t)|^2$ on the plane $z=0$}
    \label{fig:multi_coulomb_density}
\end{figure}

To simulate such dynamics, we use the ST-PP-FP method with $\tau=5\times10^{-4}$ and $h = 2^{-3}$.
% The output times are
% \[
%     t=0,\ 0.5,\ 1,\ 1.25,\ 1.5,\ 1.75,\ 2,\ 2.25,\ 2.5 .
% \]
The density $|\psi(\mathbf{x},t)|^2$ on the plane $z=0$ at different times is plotted in \cref{fig:multi_coulomb_density}. 
% The green markers indicate the positions of the three Coulomb centers. 
One observes that the initially localized packet is attracted by the singular centers, then splits and redistributes among them. The density repeatedly changes its dominant location, showing a clear transfer of electron density between different nuclei. 
% This example demonstrates that the PP-FP method can handle a strongly singular three-dimensional potential while preserving a simple Fourier pseudo-spectral implementation.

\section{Conclusion}\label{sec:5}
We proposed a Fourier pseudospectral method with potential projection (PP-FP) for the nonlinear Schr\"odinger equation (NLSE) with low-regularity and singular potentials. We first introduced the NLSE with projected potential (NLSEwPP) by replacing the potential $V$ with its Fourier projection $P_NV$ in the NLSE, and established the well-posedness and convergence of the NLSEwPP to the NLSE. The PP-FP method is obtained by applying the standard Fourier pseudospectral method to NLSEwPP on a grid that resolves the projected potential. The PP-FP method has the same computational cost as the standard FP method, and the same accuracy as the Fourier spectral method. Moreover, when combined with the time-splitting methods, it preserves the structure-preservation properties of the semidiscrete splitting scheme, which will be lost if using the Fourier spectral discretization. 

\section*{Acknowledgments} This work was partially supported by the Ministry of Education of Singapore under its AcRF Tier 1 funding grant A-8003584-00-00 (W. Bao).

\end{document}